\documentclass[numbers,webpdf]{ima-authoring-template}%

\graphicspath{{Fig/}}

\theoremstyle{thmstyletwo}%
\newtheorem{theorem}{Theorem}
\newcommand{\norm}[1]{\|#1\|}
\newcommand{\abs}[1]{|#1|}

\newcommand{\uu}{{\bf u}}
\newtheorem{remark}{Remark}%
\newtheorem{lemma}{Lemma}
\newtheorem{assumption}{Assumption}
\newtheorem{definition}{Definition}

\DeclareMathOperator{\tr}{tr}
\usepackage{algorithm}%
\usepackage{natbib}
\usepackage{moreverb}
\usepackage{algpseudocode}
\usepackage{graphicx}%
\usepackage{multirow}%
\usepackage{amsmath}%
\usepackage{amsthm}%
\usepackage{mathrsfs}%
\usepackage{xcolor}%
\usepackage{textcomp}%
\usepackage{manyfoot}%
\usepackage{booktabs}%
\usepackage{listings}%
\usepackage{float}
\floatstyle{plaintop}
\restylefloat{table}

\renewcommand{\algorithmicensure}{\textbf{Output:}}

\numberwithin{equation}{section}

\begin{document}

\DOI{DOI HERE}
\copyrightyear{2025}
\vol{00}
\pubyear{2025}
\access{Advance Access Publication Date: Day Month Year}
\appnotes{Paper}
\copyrightstatement{Published by Oxford University Press on behalf of the Institute of Mathematics and its Applications. All rights reserved.}
\firstpage{1}


\title[CholeskyQR with the randomized model]{An improved error analysis of CholeskyQR with the randomized model}

\author{Haoran Guan*
\address{\orgdiv{Department of Applied Mathematics}, \orgname{the Hong Kong Polytechnic University}, \orgaddress{\street{Hung Hom, Kowloon}, \postcode{999077}, \state{Hong Kong SAR}, \country{China}}}}
\author{Yuwei Fan
\address{\orgdiv{Theory Lab}, \orgname{Huawei Hong Kong Research Center}, \orgaddress{\street{Shatin, New Territories}, \postcode{999077}, \state{Hong Kong SAR}, \country{China}}}}

\authormark{Haoran Guan et al.}

\corresp[*]{Corresponding author: \href{email:email-id.com}{21037226R@connect.polyu.hk}}

\received{Date}{0}{Year}
\revised{Date}{0}{Year}
\accepted{Date}{0}{Year}


\abstract{This work is about an improved error analysis of CholeskyQR with the randomized model for the tall-skinny $X \in \mathbb{R}^{m\times n}$. Due to the structure of CholeskyQR, we utilize the randomized model in the first step of CholeskyQR with a weak assumption. We receive a better sufficient condition of $\kappa_{2}(X)$ and a tighter upper bound of residual for CholeskyQR2, together with a probabilistic shifted item $s$ for Shifted CholeskyQR3 based on $\norm{X}_{F}$ after improved error analysis. Numerical experiments demonstrate the effectiveness of our new theoretical results. The probabilistic $s$ for Shifted CholeskyQR3 can enhance the applicability of Shifted CholeskyQR3 while maintaining numerical stability. It is also robust enough after numerous experiments.}


\maketitle

\section{Introduction}
CholeskyQR is a very popular algorithm in both academia and industry for QR factorization. Compared with other QR algorithms, such as CGS(2), MGS(2), Householder, and TSQR \cite{ballard2011, 2011, Communication, MatrixC, Higham, Numerical}, CholeskyQR strikes a balance between efficiency, accuracy, and computational cost. Moreover, it utilizes BLAS3 operations and requires only simple reductions in a parallel environment, whereas the other algorithms require significantly more. CholeskyQR primarily addresses the case of tall-skinny matrices. If the input matrix $X \in \mathbb{R}^{m\times n}$ is full-rank, $m \ge n$, we can use CholeskyQR to perform QR factorization of $X$, see Algorithm~\ref{alg:cholqr}. Here, $Q \in \mathbb{R}^{m\times n}$ is the orthogonal factor and $R \in \mathbb{R}^{n\times n}$ is the upper-triangular factor.

\begin{algorithm}[H]
\caption{$[Q,R]=\mbox{CholeskyQR}(X)$}
\label{alg:cholqr}
\begin{algorithmic}[1]
\Require $X \in \mathbb{R}^{m\times n}.$ 
\algorithmicensure \mbox{Orthogonal factor} $Q \in \mathbb{R}^{m\times n}$, \mbox{Upper triangular factor} $R \in \mathbb{R}^{n \times n}.$ \\
 $G  = X^{\top}X,$ \\
 $R = \mbox{Cholesky}(G),$ \\
 $Q  =  XR^{-1}.$
\end{algorithmic}
\end{algorithm}%

\subsection{Some CholeskyQR-type algorithms}
Though with  many advantages, CholeskyQR is seldom used directly due to its lack of orthogonality. Therefore, a new algorithm called CholeskyQR2 was developed by applying CholeskyQR twice \cite{2014, error}, as listed in Algorithm~\ref{alg:cholqr2}. CholeskyQR2 exhibits good orthogonality and residual properties. However, it is not applicable to many ill-conditioned $X$ since the step of Cholesky factorization may break down due to rounding errors. To avoid numerical breakdown, a shifted item $s$ is added in the step of Cholesky factorization and we can get Shifted CholeskyQR (SCholeskyQR) \cite{Shifted}. Researchers put CholeskyQR2 directly after Shifted CholeskyQR to form Shifted CholeskyQR3 (SCholeskyQR3) in order to achieve numerical stability. Shifted CholeskyQR and Shifted CholeskyQR3 are detailed in Algorithm~\ref{alg:Shifted} and Algorithm~\ref{alg:Shifted3}, respectively. In addition to deterministic algorithms, some randomized algorithms for CholeskyQR have also been developed, see \cite{Randomized, Novel, Householder}.

\begin{algorithm}[H]
\caption{$[Q,R]=\mbox{CholeskyQR2}(X)$}
\label{alg:cholqr2}
\begin{algorithmic}[1]
\Require $X \in \mathbb{R}^{m\times n}.$ 
\algorithmicensure \mbox{Orthogonal factor} $Q \in \mathbb{R}^{m\times n}$, \mbox{Upper triangular factor} $R \in \mathbb{R}^{n \times n}.$ \\
 $[W,Y]=\mbox{CholeskyQR}(X),$ \\
 $[Q,Z]=\mbox{CholeskyQR}(Q),$ \\
 $R=ZY.$
\end{algorithmic}
\end{algorithm}%

\begin{algorithm}[H]
\caption{$[Q,R]=\mbox{SCholeskyQR}(X)$}
\label{alg:Shifted}
\begin{algorithmic}[1]
\Require $X \in \mathbb{R}^{m\times n}.$ 
\algorithmicensure \mbox{Orthogonal factor} $Q \in \mathbb{R}^{m\times n}$, \mbox{Upper triangular factor} $R \in \mathbb{R}^{n \times n}.$ \\
 $G  = X^{\top}X,$ \\
 take $s>0,$ \\
 $R = \mbox{Cholesky}(G+sI),$ \\
 $Q  =  XR^{-1}.$
\end{algorithmic}
\end{algorithm}%

\begin{algorithm}[H]
\caption{$[Q,R]=\mbox{SCholeskyQR3}(X)$}
\label{alg:Shifted3}
\begin{algorithmic}[1]
\Require $X \in \mathbb{R}^{m\times n}.$ 
\algorithmicensure \mbox{Orthogonal factor} $Q \in \mathbb{R}^{m\times n}$, \mbox{Upper triangular factor} $R \in \mathbb{R}^{n \times n}.$ \\
 $[W,Y]=\mbox{SCholeskyQR}(X),$ \\
 $[Q,Z]=\mbox{CholeskyQR2}(W),$ \\
 $R=ZY.$
\end{algorithmic}
\end{algorithm}%

In our recent work \cite{Columns}, we provide an improved shifted item $s$ based on $\norm{X}_{g}$ for the input $X$, a new definition based on columns. This improved $s$ can deal with more ill-conditioned $X$ while keeping numerical stability. $\norm{\cdot}_{g}$ plays an important role in rounding error analysis when $X$ has certain structure, \textit{e.g.}, $X$ is sparse. In \cite{CSparse}, we focus on Shifted CholeskyQR for sparse matrices, with an alternative $s$ based on the structure of $X$. $\norm{\cdot}_{g}$ is used in the theoretical analysis. This is the first work building connections between rounding error analysis and sparse matrices to the best of our knowledge.

\subsection{Problems and considerations}
In \cite{Columns, CSparse}, we have already made progress regarding CholeskyQR-type algorithms from the perspective of error analysis. However, we are still exploring whether we can provide more accurate error analysis for these algorithms, especially a better shifted item $s$ for Shifted CholeskyQR3. Although some improvements have been made in our previous work, there is still some distance from satisfaction.

In fact, all the works regarding CholeskyQR-type algorithms are based on the deterministic models in \cite{Higham}. However, the deterministic models of rounding error analysis will lead to an overestimation of the norms of error matrices in most of the cases, resulting in a conservative shifted items $s$ and restricted sufficient conditions for $\kappa_{2}(X)$. In floating-point arithmetic, the norms of error matrices have a very small probability of reaching the upper bounds based on the deterministic models. Randomized linear algebra has gained popularity in recent years. One of the important components of randomized linear algebra is probabilistic error analysis. Some randomized models \cite{Stochastic, New} for probabilistic error analysis of matrix multiplications have been developed, which are helpful in providing sharper upper bounds for error estimations along with the corresponding probabilities. There are some works in recent years regarding probabilistic error analysis using the randomized models, see \cite{HQR, inner, MGS} and their references. For CholeskyQR, the property of the algorithm primarily depends on the first $X^{\top}X$, especially for the tall-skinny $X \in \mathbb{R}^{m\times n}$ with $m>>n$. We are curious whether it is applicable to utilize the idea of probabilistic error analysis to provide an improved error analysis of CholeskyQR.

\subsection{Our contributions}
In this work, we use the randomized model for probabilistic error analysis in \cite{Stochastic} to provide an improved error analysis of CholeskyQR. One of the innovative points of our work is combining the idea of probabilistic error analysis with CholeskyQR-type algorithms. Different from the existing works regarding probabilistic error analysis, we use the randomized model in the first step of $X^{\top}X$ because of the structure of CholeskyQR. From the perspective of the theoretical analysis, our work introduce $\norm{\cdot}_{F}$ in the error analysis, which is also distinguished from other existing works \cite{Shifted, error, Columns, CSparse}. For CholeskyQR2, we receive a better sufficient condition of $\kappa_{2}(X)$ and a tighter error bound of residual for the input matrix $X$, which is an improvement compared to the analysis in \cite{error}. For Shifted CholeskyQR3, a probabilistic $s$ is also taken in order to improve the applicability of the algorithm while maintaining numerical stability. As far as we know, this is also the first time for the theoretical results from the error analysis based on the randomized model to be used to improve the property of a numerical algorithm. Numerical experiments demonstrate the effectiveness of our new theoretical results. The probabilistic $s$ also has better performances compared to $s$ in \cite{Columns}. The test of robustness of the probabilistic $s$ is also an innovative point of this work.

\subsection{Outline and notations}
The paper is organized as follows. In Section~\ref{sec:lit}, we present some lemmas from the existing works. We show the improved error analysis of CholeskyQR2 and Shifted CholeskyQR3 in Section~\ref{sec:ch2} and Section~\ref{sec:sc3} with the randomized model. Numerical experiments are detailed in Section~\ref{sec:experiments}, followed by conclusions in Section~\ref{sec:conclusions}.

In this work, $\norm{\cdot}_{F}$ and $\norm{\cdot}_{2}$ denote the Frobenius norm and the $2$-norm of the matrix. For $X \in \mathbb{R}^{m\times n}$ and is full-rank, $\norm{X}_{2}=\sigma_{1}(X)$ is the largest singular value and $\sigma_{min}(X)$ is the smallest singular value of $X$. $\kappa_{2}(X)=\frac{\norm{X}_{2}}{\sigma_{min}(X)}$ is the condition number of $X$. $\uu=2^{-53}$ is the unit roundoff. For the input matrix $X$, $\abs{X}$ is the matrix whose elements are all the absolute values of the elements of $X$.

\section{Some lemmas from the existing works}
\label{sec:lit}
In this section, we provide a brief review of some lemmas for theoretical analysis in this work.

\subsection{Rounding error analysis}
In the beginning of this part, we show the following classical
model for floating-point arithmetic from \cite{Higham}.
\begin{equation}
fl(a \quad \mbox{op} \quad b) = (a \quad \mbox{op} \quad b)(1+\delta), \quad \abs{\delta} \le u, \quad \mbox{op} \in \{ +, - , \times , /, \surd \}. \label{eq:model}
\end{equation}
Here, $fl(\cdot)$ denotes the computed value in floating-point arithmetic. $\delta$ is the rounding error in the operation. In the following, we show some fundamental lemmas of deterministic rounding error analysis \cite{MatrixC, Higham}, which are widely used in the error estimations of numerical linear algebra.

\begin{lemma}[Weyl's theorem of singular values]
\label{lemma 2.1}
If $A,B \in \mathbb{R}^{m\times n}$, then
\begin{equation}
\sigma_{min}(A+B) \ge \sigma_{min}(A)-\norm{B}_{2}. \nonumber
\end{equation}
\end{lemma}

\begin{lemma}[Deterministic error analysis of matrix multiplications]
\label{lemma 2.6}
For $A \in \mathbb{R}^{m\times n}, B \in \mathbb{R}^{n\times p}$, the error in computing the matrix product $AB$ in floating-point arithmetic is bounded by
\begin{equation}
\abs{AB-fl(AB)}\le \gamma_{n}\abs{A}\abs{B}. \nonumber
\end{equation}
Here, $\abs{A}$ is the matrix whose $(i,j)$ element is $\abs{a_{ij}}$ and
\begin{equation}
\gamma_n: = \frac{n{\uu}}{1-n{\uu}} \le 1.02n{\uu}. \nonumber
\end{equation}
\end{lemma}

\begin{lemma}[Deterministic error analysis of Cholesky factorization]
\label{lemma 2.8}
If Cholesky factorization applied to the symmetric positive definite $A \in \mathbb{R}^{n\times n}$ runs to completion, then the computed factor $R \in \mathbb{R}^{n\times n}$ satisfies
\begin{equation}
R^{\top}R=A+\Delta{A}, \quad \abs{\Delta A}\le \gamma_{n+1}\abs{{R}^{\top}}\abs{R}. \nonumber
\end{equation}
\end{lemma}

In the following, we present the following lemmas related to probabilistic error analysis \cite{Stochastic}. We show the randomized model of rounding errors first.

\begin{lemma}[The randomized model]
\label{lemma 2.5}
Let the computation interest generate rounding errors $\delta_{1}, \delta_{2}, \cdots$ in that order. The $\delta_{k}$ are random variables of mean zero such that $\mathbb{E}(\delta_{k}|\delta_{1}, \cdots, \delta_{k-1})$=$\mathbb{E}(\delta_{k})=0$. 
\end{lemma}

Before showing the lemma of probabilistic error analysis, we define 
\begin{align} 
P(\eta) &= 1-2\exp(-\frac{\eta^{2}(1-\uu)^{2}}{2}), \label{eq:9} \\
Q(\eta, n) &= 1-n(1-P(\eta)),  \label{eq:Q} \\
\tilde\gamma_{n} &= \exp(\eta\sqrt{n}\uu+\frac{n\uu^{2}}{1-\uu})-1. \label{eq:7}
\end{align}
Here, $\eta$ is a positive constant. From \eqref{eq:7}, we can find that when $\eta\sqrt{n}\uu$ is small and close to $0$, $\tilde\gamma_{n} \approx 1.02\eta\sqrt{n}\uu$. In the following, we present a lemma from \cite{Stochastic} regarding probabilistic error analysis in matrix multiplications, which will be used in the theoretical analysis of this work.

\begin{lemma}[Probabilistic error analysis of matrix multiplications]
\label{lemma 2.6a}
For $A \in \mathbb R^{m\times n}$ and $B \in \mathbb R^{n\times p}$, under Lemma~\ref{lemma 2.5}, the error in computing the matrix product $C=AB$ in floating-point arithmetic satisfies
\begin{equation}
\abs{AB-fl(AB)} \le \tilde{\gamma_{n}}(\eta)\abs{A}\abs{B}, \nonumber
\end{equation}
with probability at least $Q(\eta, mnp)$.
\end{lemma}

To learn more about rounding error analysis and matrix perturbation, readers can refer to \cite{Perturbation}.

\subsection{CholeskyQR2}
For CholeskyQR2 \cite{2014, error}, the main theoretical results are presented in Lemma~\ref{lemma 2.9}.

\begin{lemma}[Deterministic error analysis of CholeskyQR2]
\label{lemma 2.9}
For $X \in \mathbb{R}^{m\times n}$ and $[Q,R]=\mbox{CholeskyQR2}(X)$, when 
\begin{equation}
\delta=8\kappa_{2}(X)\sqrt{mn\uu+n(n+1)\uu} \le 1, \label{eq:15}
\end{equation}
with $mn\uu \le \frac{1}{64}$ and $n(n+1)\uu \le \frac{1}{64}$, we have
\begin{align}
\norm{Q^{\top}Q-I}_{F} &\le 6(mn\uu+n(n+1)\uu), \label{eq:13} \\
\norm{QR-X}_{F} &\le 5n^{2}\sqrt{n}\uu\norm{X}_{2}. \label{eq:14}
\end{align}
\end{lemma}

\subsection{The improved Shifted CholeskyQR3}
Shifted CholeskyQR3 can effectively handle many ill-conditioned matrices than CholeskyQR2 \cite{Shifted}. In \cite{Columns}, we provide an improved Shifted CholeskyQR3 in \cite{Columns}. We define a new matrix norm, $\norm{\cdot}_{g}$ below.

\begin{definition}
\label{def:g}
For $X=[X_{1},X_{2}, \cdots X_{n-1},X_{n}]\in R^{m\times n}$,
\begin{equation} 
\begin{split}
\norm{X}_{g}:=\max_{1 \le j \le n}\norm{X_{j}}_{2}, \nonumber
\end{split}
\end{equation}
where
\begin{equation}
\begin{split}
\norm{X_{j}}_{2}=\sqrt{x_{1,j}^{2}+x_{2,j}^{2}+……+x_{m-1,j}^{2}+x_{m,j}^{2}}. \nonumber
\end{split}
\end{equation}
\end{definition}

With $\norm{\cdot}_{g}$, we have the following theoretical results for the improved Shifted CholeskyQR3.

\begin{lemma}[The relationship between $\kappa_{2}(X)$ and $\kappa_{2}(W)$]
\label{lemma 2.13}
For $X \in \mathbb{R}^{m\times n}$ and $[W,Y]=\mbox{SCholeskyQR}(X)$, with $4.89pn^{2}\uu\kappa_{2}(X) \le 1$, $11(mn\uu+n(n+1)\uu)\norm{X}_{g}^{2} \le s \le \frac{1}{100}\norm{X}_{g}^{2}$, $mn\uu \le \frac{1}{64}$ and $n(n+1)\uu \le \frac{1}{64}$, we have
\begin{equation}
\kappa_{2}(W) \le 2\sqrt{3} \cdot \sqrt{1+t(\kappa_{2}(X))^{2}}. \label{eq:20}
\end{equation}
Here, we have $t=\frac{s}{\norm{X}_{2}^{2}}$. For $[Q,R]=\mbox{SCholeskyQR3}(X)$, when $s=11(mn\uu+n(n+1)\uu)[X]_{g}^{2}$ and $\kappa_{2}(X)$ is large enough, the sufficient condition of $\kappa_{2}(X)$ is
\begin{equation}
\kappa_{2}(X) \le \frac{1}{86p(mn\uu+n(n+1)\uu)}. \label{eq:23}
\end{equation}
Here, we let $p=\frac{\norm{X}_{g}}{\norm{X}_{2}}$, $\frac{1}{\sqrt{n}} \le p \le 1$.
\end{lemma}

\begin{lemma}[Deterministic error analysis of the improved Shifted CholeskyQR3]
\label{lemma 2.14}
For $X \in \mathbb{R}^{m\times n}$ and $[Q,R]=\mbox{SCholeskyQR3}(X)$, when $s=11(mn\uu+n(n+1)\uu)\norm{X}_{g}^{2}$, with $mn\uu \le \frac{1}{64}$, $n(n+1)\uu \le \frac{1}{64}$ and \eqref{eq:23}, we have
\begin{align}
\norm{Q^{\top}Q-I}_{F} &\le 6(mn\uu+n(n+1)\uu), \label{eq:21} \\
\norm{QR-X}_{F} &\le (6.57p+4.81)n^{2}\uu\norm{X}_{2}. \label{eq:22}
\end{align}
\end{lemma}

\section{An improved error analysis of CholeskyQR2 with the randomized model} 
\label{sec:ch2}
In this section, we aim to utilize the randomized model to conduct an improved error analysis of CholeskyQR2. Before presenting the theoretical results of this work, we make an assumption first.

\begin{assumption} \label{assumption 1}
The first step of CholeskyQR for the input matrix $X \in \mathbb{R}^{m\times n}$, $X^{\top}X$, follows Lemma~\ref{lemma 2.5} and Lemma~\ref{lemma 2.6a}.
\end{assumption}

Actually, it is hard for us to do probabilistic error analysis for the whole CholeskyQR because the algorithm consists of several independent steps of computation one by one. Assumption~\ref{assumption 1} is a weak assumption which only focuses on the first step of $X^{\top}X$. However, for the tall-skinny $X \in \mathbb{R}^{m\times n}$ with $m>>n$ and the new technique of error analysis, it is enough for us provide better results compared to those in the existing works.

\subsection{General settings}
\label{sec:g1}
In the beginning, we present CholeskyQR2 step by step, accompanied by the corresponding error matrices below.
\begin{align}
G-X^{\top}X &=E_{A}, \label{eq:35} \\
Y^{\top}Y-G &= E_{B}, \label{eq:36} \\
WY=X &+ E_{WY}, \label{eq:37} \\
C-W^{\top}W &= E_{1}, \nonumber \\
Z^{\top}Z-C &= E_{2}, \label{eq:39} \\
QZ-W &= E_{3}, \label{eq:310} \\
ZY-R &= E_{4}. \label{eq:311}
\end{align}

For the input matrix $X \in \mathbb{R}^{m\times n}$, we provide some general settings for all the algorithms in this chapter below.
\begin{align}
\max(\eta \sqrt{m}n\uu, mn\uu) \le \frac{1}{64}, \label{eq:31} \\
n(n+1)\uu \le \frac{1}{64}. \label{eq:32}
\end{align}
Here, $\eta$ occurs in \eqref{eq:9} and \eqref{eq:Q}. For CholeskyQR2, when \eqref{eq:31} and \eqref{eq:32} are satisfied, if we want to have $Q(\eta, mn^{2})$ to be positive, we can choose $\eta$ as a positive constant not exceeding $10$ in numerical experiments. Moreover, we let 
\begin{equation} 
j=\frac{\norm{X}_{F}}{\norm{X}_{2}}. \label{eq:kxf2}
\end{equation}
Here, $1 \le j \le \sqrt{n}$. 

\subsection{The improved error analysis of CholeskyQR2 based on the randomized model}
In this section, we present some theoretical results related to the improved error analysis of CholeskyQR2 with the randomized model partially.

\begin{theorem}[The improved error analysis of CholeskyQR2 with the randomized model]
\label{theorem 3.1}
Under Assumption~\ref{assumption 1}, for $X \in \mathbb{R}^{m\times n}$ and $[Q,R]=\mbox{CholeskyQR2}(X)$, with \eqref{eq:31}, \eqref{eq:32} and 
\begin{equation}
\delta_{1}=8j\kappa_{2}(X)\sqrt{\eta\sqrt{m}\uu+(n+1)\uu} \le 1, \label{eq:deltac2}
\end{equation}
we have
\begin{align}
\norm{Q^{\top}Q-I}_{F} &\le 6(mn\uu+n(n+1)\uu), \label{eq:33} \\
\norm{QR-X}_{F} &\le (2.30j+1.21\sqrt{n}) \cdot n\sqrt{n}\uu\norm{X}_{2}, \label{eq:34}
\end{align}
with probability at least $Q(\eta, mn^{2})$. $j$ is defined in \eqref{eq:kxf2}. 
\end{theorem}

The comparison between the theoretical results of CholeskyQR2 under different types of analysis are shown below in Table~\ref{tab:ComparisonC}. We refer to the theoretical results in Lemma~\ref{lemma 2.9} from \cite{error} as ‘\textit{Deterministic}’ and our new results in Theorem~\ref{theorem 3.1} as ‘\textit{Probabilistic}’.

\begin{table}
\caption{Comparison of the theoretical results of CholeskyQR2 between different analysis for $X \in \mathbb{R}^{m\times n}$}
\centering
\begin{tabular}{||c c c||}
\hline
$\mbox{Type of analysis}$ & $\mbox{Condition of $\kappa_{2}(X)$}$ & $\norm{QR-X}_{F}$ \\
\hline
$\mbox{Deterministic}$ & $\frac{1}{8\sqrt{mn\uu+n(n+1)\uu}}$ & $5n^{2}\sqrt{n}\uu\norm{X}_{2}$ \\
\hline
$\mbox{Probabilistic}$ & $\frac{1}{8j\sqrt{\eta\sqrt{m}\uu+(n+1)\uu}}$ & $(2.30j+1.21\sqrt{n}) \cdot n\sqrt{n}\uu\norm{X}_{2}$ \\
\hline
\end{tabular}
\label{tab:ComparisonC}
\end{table}

\subsection{Lemmas for proving Theorem~\ref{theorem 3.1}}
Before proving Theorem~\ref{theorem 3.1}, we present some lemmas related to it. The analytic steps of some lemmas in this chapter are similar to those in \cite{Shifted, error, Columns, CSparse}. However, we utilize the randomized model and $\norm{\cdot}_{F}$, allowing us to obtain sharper upper bounds with high probabilities, which are fundamentally different from these existing works. 

\begin{lemma}
\label{lemma 3.1}
For $E_{A}$ and $E_{B}$ in \eqref{eq:35} and \eqref{eq:36}, we have
\begin{align}
\norm{E_{A}}_{2} &\le 1.1\eta \sqrt{m}\uu\norm{X}_{F}^{2}, \label{eq:eac} \\
\norm{E_{B}}_{2} &\le 1.1(n+1)\uu\norm{X}_{F}^{2}, \label{eq:ebc}
\end{align}
with probability at least $Q(\eta, mn^{2})$.
\end{lemma}
\begin{proof}
Regarding $\norm{E_{A}}_{2}$, with Lemma~\ref{lemma 2.6a} and \eqref{eq:35}, we can have
\begin{equation} \label{eq:312}
\begin{split}
\abs{E_{A}} &= \abs{G-X^{\top}X} \\ &\le \tilde\gamma_{m}(\eta)\abs{X^{\top}}\abs{X} \\ &\le 1.1\eta\sqrt{m}\uu \cdot \abs{X^{\top}}\abs{X}, 
\end{split}
\end{equation}
with probability at least $Q(\eta, mn^{2})$. Therefore, we can bound $\norm{E_{A}}_{2}$ as 
\begin{equation}
\begin{split}
\norm{E_{A}}_{2} &\le \norm{\abs{E_{A}}}_{F} \\ &\le \tilde\gamma_{m}(\eta) \cdot \norm{X}_{F}^{2} \\ &\le 1.1\eta\sqrt{m}\uu\norm{X}_{F}^{2}, \nonumber
\end{split}
\end{equation}
with probability at least $Q(\eta, mn^{2})$. \eqref{eq:eac} is proved.

Regarding $\norm{E_{B}}_{2}$, with Lemma~\ref{lemma 2.8}, \eqref{eq:35} and \eqref{eq:36}, we can have
\begin{equation} \label{eq:314}
\begin{split}
\norm{E_{B}}_{2} &\le \norm{\abs{E_{B}}}_{F} \\ &\le \gamma_{n+1} \cdot \norm{Y}_{F}^{2} \\ &\le \gamma_{n+1} \cdot \tr(Y^{\top}Y) \\ &\le \gamma_{n+1} \cdot \tr(X^{\top}X+E_{A}+E_{B}) \\ &\le \gamma_{n+1} \cdot (\norm{X}_{F}^{2}+n\norm{E_{A}}_{2}+n\norm{E_{B}}_{2}), 
\end{split}
\end{equation}
with probability at least $Q(\eta, mn^{2})$. With \eqref{eq:31}, \eqref{eq:32}, \eqref{eq:eac} and \eqref{eq:314}, we can get
\begin{equation}
\begin{split}
\norm{E_{B}}_{2} &\le \frac{\gamma_{n+1}(1+1.1\tilde\gamma_{m}(\eta) \cdot n)}{1-\gamma_{n+1} \cdot n} \cdot \norm{X}_{F}^{2} \nonumber \\
&\le \frac{1.02(n+1)\uu \cdot (1+1.1\eta \sqrt{m}\uu)}{1-1.02(n+1)\uu \cdot n} \cdot \norm{X}_{F}^{2} \nonumber \\
&\le \frac{1.02(n+1)\uu \cdot (1+1.1\cdot\frac{1}{64})}{1-\frac{1.02}{64}} \cdot \norm{X}_{F}^{2} \nonumber \\ 
&\le 1.1(n+1)\uu\norm{X}_{F}^{2}, \nonumber
\end{split}
\end{equation}
with probability at least $Q(\eta, mn^{2})$. \eqref{eq:ebc} is proved. Therefore, Lemma~\ref{lemma 3.1} holds.
\end{proof}

\begin{lemma}
\label{lemma 3.2}
For $Y^{-1}$ and $XY^{-1}$ in \eqref{eq:37}, we have
\begin{align}
\norm{Y^{-1}}_{2} &\le \frac{1.1}{\sigma_{min}(X)}, \label{eq:322} \\
\norm{XY^{-1}}_{2} &\le 1.1, \label{eq:323}
\end{align}
\end{lemma}
\begin{proof}
with probability at least $Q(\eta, mn^{2})$. The idea to prove Lemma~\ref{lemma 3.2} is the same as that in \cite{error}. Based on Lemma~\ref{lemma 2.1}, \eqref{eq:35} and \eqref{eq:36}, we can have
\begin{equation}
(\sigma_{min}(Y))^{2} \ge (\sigma_{min}(X))^{2}-(\norm{E_{A}}_{2}+\norm{E_{B}}_{2}). \label{eq:324}
\end{equation}
Therefore, based on \eqref{eq:kxf2}, \eqref{eq:deltac2} and \eqref{eq:ebc}, we can have
\begin{equation} \label{eq:325}
\begin{split}
\norm{E_{A}}_{2}+\norm{E_{B}}_{2} &\le \frac{1.1}{64}(\sigma_{min}(X))^{2} \\ &\le (1-\frac{1}{1.1^{2}})(\sigma_{min}(X))^{2}, 
\end{split}
\end{equation}
with probability at least $Q(\eta, mn^{2})$. We combine \eqref{eq:324} with \eqref{eq:325} and we can have
\begin{equation}
\frac{1}{1.1^{2}} \cdot (\sigma_{min}(X))^{2} \le (\sigma_{min}(Y))^{2}, \label{eq:326} 
\end{equation}
with probability at least $Q(\eta, mn^{2})$. Therefore, we can easily get \eqref{eq:322}. The same as the step \cite{error}, we can have \eqref{eq:323}. Lemma~\ref{lemma 3.2} holds.
\end{proof}

\begin{lemma}
\label{lemma 3.3}
For $E_{WY}$ in \eqref{eq:37}, we have
\begin{equation} \label{eq:331}
\norm{E_{WY}}_{2} \le 1.05n\sqrt{n}\uu \cdot \norm{W}_{2}\norm{X}_{F}, 
\end{equation}
with probability at least $Q(\eta, mn^{2})$. 
\end{lemma}
\begin{proof}
With Lemma~\ref{lemma 2.6} and \eqref{eq:37}, we can have
\begin{equation} \label{eq:ewyc2}
\begin{split}
\norm{E_{WY}}_{2} &\le \gamma_{n} \cdot \norm{W}_{F}\norm{Y}_{F} \\ &\le 1.02n\sqrt{n}\uu \cdot \norm{W}_{2}\norm{Y}_{F}.
\end{split}
\end{equation}
Actually, if we do not use the randomized model, we can easily receive the deterministic upper bound of $\norm{E_{A}}_{2}$ as $1.1m\uu\norm{X}_{F}^{2}$. Therefore, similar to the steps to get \eqref{eq:ebc}, we can get the deterministic upper bounds of both $\norm{E_{B}}_{2}$ as $1.1(n+1)\uu\norm{X}_{F}^{2}$. Based on \eqref{eq:35}, \eqref{eq:36}, \eqref{eq:31}, \eqref{eq:32}, \eqref{eq:ebc} and \eqref{eq:314}, we can have
\begin{equation} \label{eq:329} 
\begin{split}
\norm{Y}_{F}^{2} &=\tr(Y^{\top}Y) \\ &=\tr(X^{\top}X+E_{A}+E_{B}) \\ &\le \norm{X}_{F}^{2}+n(\norm{E_{A}}_{2}+\norm{E_{B}}_{2}) \\ &\le 1.04\norm{X}_{F}^{2}, 
\end{split}
\end{equation}
with probability at least $Q(\eta, mn^{2})$. Based on \eqref{eq:329}, we can have 
\begin{equation}
\norm{Y}_{F} \le 1.02\norm{X}_{F}, \label{eq:330}
\end{equation}
with probability at least $Q(\eta, mn^{2})$. Therefore, we put \eqref{eq:330} into \eqref{eq:ewyc2} and we can get \eqref{eq:331}. Lemma~\ref{lemma 3.3} holds.
\end{proof}

\begin{lemma}
\label{lemma 3.4}
For $W$ in \eqref{eq:37}, we have
\begin{equation} \label{eq:341b}
\norm{W}_{2} \le 1.13, 
\end{equation}
with probability at least $Q(\eta, mn^{2})$.
\end{lemma}
\begin{proof}
Based on \eqref{eq:37}, we can have
\begin{equation} \label{eq:wfwc}
\begin{split}
\norm{W}_{2} &\le \norm{XY^{-1}}_{2}+\norm{E_{WY}Y^{-1}}_{2} \\ &\le \norm{XY^{-1}}_{2}+\norm{E_{WY}}_{2}\norm{Y^{-1}}_{2}.
\end{split}
\end{equation}
With \eqref{eq:32}, \eqref{eq:kxf2}, \eqref{eq:322} and \eqref{eq:331}, we can have
\begin{equation} \label{eq:ewyy-1c}
\begin{split}
\norm{E_{WY}}_{2}\norm{Y^{-1}}_{2} &\le \frac{1.16n\sqrt{n}\uu \cdot \norm{W}_{2}\norm{X}_{F}}{\sigma_{min}(X)} \\ &\le \frac{1.16jn\sqrt{n}\uu \cdot \norm{X}_{2}}{8j\sqrt{\eta\sqrt{m}\uu+(n+1)\uu} \cdot \norm{X}_{2}} \cdot \norm{W}_{2} \\ &\le \frac{1.16}{8} \cdot n\sqrt{\uu} \cdot \norm{W}_{2} \\ &\le 0.02\norm{W}_{2},  
\end{split}
\end{equation}
with probability at least $Q(\eta, mn^{2})$. Therefore, we put \eqref{eq:323} and \eqref{eq:ewyy-1c} into \eqref{eq:wfwc} and we can have
\begin{equation} \label{eq:w2vc}
\norm{W}_{2} \le 1.1+0.02\norm{W}_{2}, 
\end{equation}
with probability at least $Q(\eta, mn^{2})$. With \eqref{eq:w2vc}, we can have \eqref{eq:341b}. Lemma~\ref{lemma 3.4} holds.
\end{proof}
\begin{remark}
In Lemma~\ref{lemma 3.3} and Lemma~\ref{lemma 3.4}, we use a distinguished idea to do theoretical analysis regarding $W$ and $E_{WY}$ in \eqref{eq:37}. In \cite{Shifted, error, Columns, CSparse}, $E_{WY}$ is estimated through the idea of solving a linear system for each row of the input $X$. This idea is applicable when we take the deterministic error analysis. However, when we consider the randomized model, the corresponding steps in these works leads to a component of the least probability as $(Q(\eta, \frac{n(n+1)}{2}))^{m}$ if the input matrix is $X \in \mathbb{R}^{m\times n}$. It is not acceptable if we want to have a high least probability. Our new analysis in this work addresses this problem in a more direct way of error estimation, receiving a proper least probability $Q(\eta, mn^{2})$.
\end{remark}

\subsection{Proof of Theorem~\ref{theorem 3.1}}
With Lemma~\ref{lemma 3.1}-Lemma~\ref{lemma 3.4}, we begin to prove Theorem~\ref{theorem 3.1}.

\begin{proof}
The proof of Theorem~\ref{theorem 3.1} is divided into two parts, orthogonality and residual.

\subsubsection{The upper bound of orthogonality}
First, we consider the orthogonality. Based on \eqref{eq:35}, \eqref{eq:36} and \eqref{eq:37}, it is easy to get
\begin{equation}
\begin{split}
W^{\top}W &= Y^{-\top}(X+E_{WY})^{\top}(X+E_{WY})Y^{-1} \nonumber \\ 
&= Y^{-\top}X^{\top}XY^{-1}+Y^{-\top}X^{\top}E_{WY}Y^{-1} \nonumber \\ 
&+ Y^{-\top}E_{WY}^{\top}XY^{-1}+Y^{-\top}E_{WY}^{\top}E_{WY}Y^{-1} \nonumber \\
&= I-Y^{-\top}(E_{A}+E_{B})Y^{-1}+(XY^{-1})^{\top}E_{WY}Y^{-1} \nonumber \\ 
&+ Y^{-\top}E_{WY}^{\top}(XR^{-1})+Y^{-\top}E_{WY}^{\top}E_{WY}Y^{-1}. \nonumber
\end{split}
\end{equation}
Therefore, we can have
\begin{equation} \label{eq:336}
\begin{split}
\norm{W^{\top}W-I}_{2} &\le \norm{Y^{-1}}_{2}^{2}(\norm{E_{A}}_{2}+\norm{E_{B}}_{2})+2\norm{Y^{-1}}_{2}\norm{XY^{-1}}_{2}\norm{E_{WY}}_{2} \\ &+ \norm{Y^{-1}}_{2}^{2}\norm{E_{WY}}_{2}^{2}.
\end{split}
\end{equation}
Based on \eqref{eq:kxf2}, \eqref{eq:deltac2}, \eqref{eq:eac}, \eqref{eq:ebc} and \eqref{eq:322}, we can have
\begin{equation} \label{eq:337} 
\begin{split}
\norm{Y^{-1}}_{2}^{2}(\norm{E_{A}}_{2}+\norm{E_{B}}_{2}) &\le \frac{1.21 \cdot (1.1j^{2} \cdot (\eta\sqrt{m}\uu+(n+1)\uu)\norm{X}_{2}^{2})}{(\sigma_{min}(X))^{2}} \\ &\le \frac{1.34}{64}\delta_{1}^{2}, 
\end{split}
\end{equation}
with probability at least $Q(\eta, mn^{2})$. Based on \eqref{eq:kxf2}, \eqref{eq:deltac2}, \eqref{eq:322}, \eqref{eq:323}, \eqref{eq:331} and \eqref{eq:341b}, we can have
\begin{equation} \label{eq:338} 
\begin{split}
2\norm{Y^{-1}}_{2}\norm{XY^{-1}}_{2}\norm{E_{WY}}_{2} &\le 2 \cdot \frac{1.1}{\sigma_{min}(X)} \cdot 1.1 \cdot (1.05jn\sqrt{n}\uu \cdot \norm{X}_{2} \cdot 1.13) \\ &\le \frac{3}{64}\delta_{1},
\end{split}
\end{equation}
with probability at least $Q(\eta, mn^{2})$. With \eqref{eq:kxf2}, \eqref{eq:deltac2}, \eqref{eq:322}, \eqref{eq:331} and \eqref{eq:341b}, we can have
\begin{equation} \label{eq:339}
\begin{split}
\norm{Y^{-1}}_{2}^{2}\norm{E_{WY}}_{2}^{2} &\le \frac{1.21}{(\sigma_{min}(X))^{2}} \cdot (1.05jn\sqrt{n}\uu \cdot \norm{X}_{2} \cdot 1.13)^{2} \\ &\le \frac{1}{2048}\delta_{1}^{2}, 
\end{split}
\end{equation}
with probability at least $Q(\eta, mn^{2})$. Therefore, we put \eqref{eq:337}-\eqref{eq:339} into \eqref{eq:336} and with \eqref{eq:deltac2}, we can have
\begin{equation}
\norm{W^{\top}W-I}_{2} \le \frac{5}{64}, \label{eq:340}
\end{equation}
with probability at least $Q(\eta, mn^{2})$. With \eqref{eq:340}, it is easy to have
\begin{align}
\norm{W}_{2} \le \frac{\sqrt{69}}{8}, \label{eq:341} \\
\sigma_{min}(W) \ge \frac{\sqrt{59}}{8}, \label{eq:342}
\end{align}
with probability at least $Q(\eta, mn^{2})$. With \eqref{eq:341} and \eqref{eq:342}, we can get
\begin{equation}
\kappa_{2}(W) \le \sqrt{\frac{69}{59}}, \label{eq:343}
\end{equation}
with probability at least $Q(\eta, mn^{2})$. Based on \eqref{eq:31}, \eqref{eq:32} and \eqref{eq:343}, we can get
\begin{equation}
j_{2}=8\kappa_{2}(W)\sqrt{mn\uu+n(n+1)\uu} \le 1, \label{eq:k2}
\end{equation}
with probability at least $Q(\eta, mn^{2})$. With \eqref{eq:k2} and similar to the result in \cite{error}, we can have
\begin{equation}
\begin{split}
\norm{Q^{\top}Q-I}_{F} &\le \frac{5}{64}j_{2}^{2} \nonumber \\
&\le 6(mn\uu+n(n+1)\uu), \nonumber
\end{split}
\end{equation}
with probability at least $Q(\eta, mn^{2})$. \eqref{eq:33} holds.

\subsubsection{The upper bound of residual}
Regarding the residual, according to \eqref{eq:310} and \eqref{eq:311}, we can have
\begin{equation}
\begin{split}
QR-X &= (W+E_{3})Z^{-1}(ZY-E_{4})-X \nonumber \\ &= (W+E_{3})Y-(W+E_{3})Z^{-1}E_{4}-X \nonumber \\ &= WY-X+E_{3}Y-QE_{4}. \nonumber
\end{split}
\end{equation}
Therefore, it is easy to have
\begin{equation}
\norm{QR-X}_{F} \le \norm{WY-X}_{F}+\norm{E_{3}}_{F}\norm{Y}_{2}+\norm{Q}_{2}\norm{E_{4}}_{F}. \label{eq:344}
\end{equation}
Based on Lemma~\ref{lemma 2.6}, \eqref{eq:kxf2}, \eqref{eq:330} and \eqref{eq:341}, we can bound $\norm{E_{WY}}_{F}=\norm{WY-X}_{F}$ as
\begin{equation} \label{eq:345}
\begin{split}
\norm{WY-X}_{F} &= \norm{E_{WY}}_{F} \\ &\le \gamma_{n} \cdot \norm{W}_{F}\norm{Y}_{F} \\ &\le 1.02n\uu \cdot \frac{\sqrt{69n}}{8} \cdot 1.02j\norm{X}_{2} \\ &\le 1.09jn\sqrt{n}\uu\norm{X}_{2}, 
\end{split}
\end{equation}
with probability at least $Q(\eta, mn^{2})$. Similar to the steps in \cite{error}, we can get
\begin{equation}
\norm{Y}_{2} \le 1.02\norm{X}_{2}. \label{eq:v2}
\end{equation}
For $Q$ in \eqref{eq:310}, with \eqref{eq:31}, \eqref{eq:32} and \eqref{eq:33}, we can have
\begin{equation}
\norm{Q}_{2} \le 1.09, \label{eq:346}
\end{equation}
with probability at least $Q(\eta, mn^{2})$. For $Z$ in \eqref{eq:39}, Similar to the steps to get \eqref{eq:330} and with \eqref{eq:341}, we can get
\begin{equation} \label{eq:349}
\begin{split}
\norm{Z}_{F} &\le 1.02\norm{W}_{F} \\ &\le 1.06\sqrt{n}, 
\end{split}
\end{equation}
with probability at least $Q(\eta, mn^{2})$. Therefore, for $E_{3}$ in \eqref{eq:310}, with Lemma~\ref{lemma 2.6}, \eqref{eq:346} and \eqref{eq:349}, we can bound $\norm{E_{3}}_{F}$ as
\begin{equation} \label{eq:347}
\begin{split}
\norm{E_{3}}_{F} &\le \gamma_{n} \cdot \norm{Q}_{F}\norm{Z}_{F} \\ &\le 1.02n\uu \cdot 1.06\sqrt{n} \cdot 1.09\sqrt{n} \\ &\le 1.18n^{2}\uu,
\end{split}
\end{equation}
with probability at least $Q(\eta, mn^{2})$. For $E_{4}$ in \eqref{eq:311}, with Lemma~\ref{lemma 2.6}, \eqref{eq:330} and \eqref{eq:349}, we can bound $\norm{E_{4}}_{F}$ as
\begin{equation} \label{eq:350}
\begin{split}
\norm{E_{4}}_{F} &\le \gamma_{n} \cdot \norm{Z}_{F}\norm{Y}_{F} \\ &\le 1.02n\uu \cdot 1.06\sqrt{n} \cdot 1.02j\norm{X}_{2} \\ & \le 1.11jn\sqrt{n}\uu\norm{X}_{2},
\end{split}
\end{equation}
with probability at least $Q(\eta, mn^{2})$. We put \eqref{eq:345}-\eqref{eq:346}, \eqref{eq:347} and \eqref{eq:350} into \eqref{eq:344} and we can have \eqref{eq:34} with probability at least $(Q(\eta, mn^{2})$. Therefore, Theorem~\ref{theorem 3.1} holds.
\end{proof}
\begin{remark}
In this part, we utilize the randomized model to provide an improved error analysis of CholeskyQR2 compared to that in \cite{error} under a weak assumption. A tighter upper bound of residual and a better sufficient condition for $\kappa_{2}(X)$ are also received, see the comparison in Table~\ref{tab:ComparisonC}. The progress can be reflected by some numerical experiments in Section~\ref{sec:experiments}.
\end{remark}

\section{An improved error analysis for Shifted CholeskyQR3 with a probabilistic $s$}
\label{sec:sc3}
We have already provided an improved error analysis for CholeskyQR2 in Section~\ref{sec:ch2}. In this part, we choose Shifted CholeskyQR3 to demonstrate the effectiveness of the randomized model and the new method of analysis. A new probabilistic shifted item $s$ is provided to improve the property of Shifted CholeskyQR3 algorithm. We have not seen similar works regarding utilizing the results of probabilistic error analysis to improve the performance of an algorithm before.

\subsection{General settings and algorithms}
\label{sec:g2}
In the beginning, we write Shifted CholeskyQR3 with error matrices step by step below. It is the same as that in \cite{Columns, CSparse}.
\begin{align}
G-X^{\top}X &= E_{A}, \label{eq:47} \\
Y^{\top}Y=G &+ sI+E_{B}, \label{eq:48} \\
WY=X &+ E_{WY}, \label{eq:49} \\
C-W^{\top}W &= E_{1}, \nonumber \\
D^{\top}D-C &= E_{2}, \label{eq:446} \\
VD-W &= E_{3}, \label{eq:447} \\
DY-N &= E_{4}, \label{eq:448} \\
B-V^{\top}V &= E_{5}, \nonumber \\
J^{\top}J-B &= E_{6}, \label{eq:450} \\
QJ-V &= E_{7}, \label{eq:451} \\
JN-R &= E_{8}. \label{eq:452}
\end{align}

For Shifted CholeskyQR3, \eqref{eq:31}-\eqref{eq:kxf2} still hold. We present more general settings below.
\begin{align}
\kappa_{2}(X) &\le \frac{1}{4.78jn\sqrt{n}\uu}, \label{eq:43} \\
11(K_{1}+(n+1)\uu)\norm{X}_{F}^{2} &\le s \le \frac{1}{100n}\norm{X}_{F}^{2}, \label{eq:44} \\
K_{2}+(n+1)\uu &\le \frac{1}{1100n}. \label{eq:k2u}
\end{align}
Here, we define $K_{1}=\min(\eta\sqrt{m}\uu,m\uu)$ and $K_{2}=\max(\eta\sqrt{m}\uu, m\uu)$.

\subsection{Probabilistic error analysis of Shifted CholeskyQR3}
In this section, we present the improved error analysis of Shifted CholeskyQR3, especially a new probabilistic $s$ based on the randomized model and $\norm{X}_{F}$ for the input matrix $X \in \mathbb{R}^{m\times n}$. 

The same as the corresponding steps in \cite{Columns, CSparse}, we separate the calculation of the $R$-factor of Shifted CholeskyQR3 into \eqref{eq:448} and \eqref{eq:452}. In the following, we show some theoretical results of Shifted CholeskyQR3 based on the new analysis.

\begin{theorem}[The relationship between $\kappa_{2}(X)$ and $\kappa_{2}(W)$]
\label{theorem 4.5}
Under Assumption~\ref{assumption 1}, for $X \in \mathbb{R}^{m\times n}$ and $[W,Y]=\mbox{SCholeskyQR}(X)$, with \eqref{eq:31}, \eqref{eq:32} and \eqref{eq:43}-\eqref{eq:k2u}, we have
\begin{equation}
\kappa_{2}(W) \le 2.88\sqrt{1+t(\kappa_{2}(X))^{2}}, \label{eq:436c}
\end{equation}
with probability at least $Q(\eta, mn^{2})$. Here, $t=\frac{s}{\norm{X}_{2}^{2}}$. When $[Q,R]=\mbox{SCholeskyQR3}(X)$, when $s=11(\eta \sqrt{m}\uu+(n+1)\uu)\norm{X}_{F}^{2}$ and $\kappa_{2}(X)$ is large enough, a sufficient condition for $\kappa_{2}(X)$ is
\begin{equation} \label{eq:2s4}
\begin{split}
\kappa_{2}(X) \le \frac{1}{77j\sqrt{(\eta\sqrt{m}\uu+(n+1)\uu)} \cdot \sqrt{(mn\uu+n(n+1)\uu)}}. 
\end{split}
\end{equation}
Here, $j$ is defined in \eqref{eq:kxf2}.
\end{theorem}

\begin{theorem}[The improved error analysis of Shifted CholeskyQR3 with the randomized model]
\label{theorem 4.6}
Under Assumption~\ref{assumption 1}, for $X \in \mathbb{R}^{m\times n}$ and $[Q,R]=\mbox{SCholeskyQR3}(X)$, when $s=11\eta(\sqrt{m}\uu+(n+1)\uu)\norm{X}_{F}^{2}$ and $\kappa_{2}(X)$ is large enough, with \eqref{eq:31}, \eqref{eq:32}, \eqref{eq:k2u} and \eqref{eq:2s4}, we have 
\begin{align}
\norm{Q^{\top}Q-I}_{F} &\le 6(mn\uu+n(n+1)\uu), \label{eq:443} \\
\norm{QR-X}_{F} &\le (5.08j+3.46\sqrt{n}) \cdot n\sqrt{n}\uu\norm{X}_{2}, \label{eq:444}
\end{align}
with probability at least $Q(\eta, mn^{2})$. Here, $j$ is defined in \eqref{eq:kxf2}.
\end{theorem}

Similar to Table~\ref{tab:ComparisonC}, comparisons of the theoretical results between different types of analysis for Shifted CholeskyQR3 are shown in Table~\ref{tab:Comparisonr3} and Table~\ref{tab:Comparisonr4}. We refer to the theoretical results in Lemma~\ref{lemma 2.14} from \cite{Columns} as ‘\textit{Deterministic}’ and the theoretical results in this work as ‘\textit{Probabilistic}’.

\begin{table}
\caption{Comparison of $\kappa_{2}(X)$ for Shifted CholeskyQR3 between different types of analysis for $X \in \mathbb{R}^{m\times n}$}
\centering
\begin{tabular}{||c c c||}
\hline
$\mbox{Type of analysis}$ & $\mbox{Sufficient condition of $\kappa_{2}(X)$}$ & $\mbox{Upper bound of $\kappa_{2}(X)$}$ \\
\hline
$\mbox{Deterministic}$ & $\frac{1}{86p(mn\uu+n(n+1)\uu)}$ & $\frac{1}{4.89pn^{2}\uu}$ \\
\hline
$\mbox{Probabilistic}$ & $\frac{1}{77j\sqrt{(\eta\sqrt{m}\uu+(n+1)\uu)} \cdot \sqrt{(mn\uu+n(n+1)\uu)}}$ & $\frac{1}{4.78jn\sqrt{n}\uu}$ \\
\hline
\end{tabular}
\label{tab:Comparisonr3}
\end{table}

\begin{table}
\caption{Comparison of the upper bound of residual for Shifted CholeskyQR3 between different types of analysis for $X \in \mathbb{R}^{m\times n}$}
\centering
\begin{tabular}{||c c||}
\hline
$\mbox{Type of analysis}$ & $\norm{QR-X}_{F}$ \\
\hline
$\mbox{Deterministic}$ & $(6.57p+4.87)n^{2}\uu\norm{X}_{2}$ \\
\hline
$\mbox{Probabilistic}$ & $(5.08j+3.46\sqrt{n}) \cdot n\sqrt{n}\uu\norm{X}_{2}$ \\
\hline
\end{tabular}
\label{tab:Comparisonr4}
\end{table}

\subsection{Lemmas for proving Theorem~\ref{theorem 4.5} and Theorem~\ref{theorem 4.6}}
To prove Theorem~\ref{theorem 4.5} and Theorem~\ref{theorem 4.6}, we present the following lemmas.

\begin{lemma}
\label{lemma 4.1}
For $E_{A}$ and $E_{B}$ in \eqref{eq:47} and \eqref{eq:48}, when \eqref{eq:44} is satisfied, we have
\begin{align}
\norm{E_{A}}_{2} &\le 1.1\eta \sqrt{m}\uu\norm{X}_{F}^{2}, \label{eq:ecc} \\
\norm{E_{B}}_{2} &\le 1.1(n+1)\uu\norm{X}_{F}^{2}, \label{eq:edc}
\end{align}
with probability at least $Q(\eta, mn^{2})$.
\end{lemma}
\begin{proof}
When estimating $\norm{E_{A}}_{2}$, the same as Lemma~\ref{lemma 3.1}, we can get \eqref{eq:ecc} with probability at least $Q(\eta, mn^{2})$.

Regarding $\norm{E_{B}}_{2}$, with Lemma~\ref{lemma 2.8}, \eqref{eq:47} and \eqref{eq:48}, we can get
\begin{equation} \label{eq:415}
\begin{split}
\norm{E_{B}}_{2} &\le \norm{\abs{E_{B}}}_{F} \\ &\le \gamma_{n+1} \cdot \norm{Y}_{F}^{2} \\ &\le \gamma_{n+1} \cdot \tr(Y^{\top}Y) \\ &\le \gamma_{n+1} \cdot \tr(X^{\top}X+sI+E_{A}+E_{B}) \\ &\le \gamma_{n+1} \cdot (\norm{X}_{F}^{2}+sn+n\norm{E_{A}}_{2}+n\norm{E_{B}}_{2}), 
\end{split}
\end{equation}
with probability at least $Q(\eta, mn^{2})$. Here, we define $t_{1}=\frac{s}{\norm{X}_{F}^{2}}$. With \eqref{eq:31}, \eqref{eq:32}, \eqref{eq:44}, \eqref{eq:ecc} and \eqref{eq:415}, we can bound $\norm{E_{B}}_{2}$ as
\begin{equation}
\begin{split}
\norm{E_{B}}_{2} &\le \frac{\gamma_{n+1} \cdot (1+\tilde\gamma_{m}(\eta) \cdot n+nt_{1})}{1-\gamma_{n+1} \cdot n}\norm{X}_{F}^{2} \nonumber \\
&\le \frac{1.02(n+1)\uu \cdot (1+1.1\eta\sqrt{m}\uu \cdot n+nt_{1})}{1-1.02(n+1)n\uu}\norm{X}_{F}^{2} \nonumber \\
&\le \frac{1.02(n+1)\uu \cdot (1+1.1\cdot\frac{1}{64}+0.01)}{1-\frac{1.02}{64}}\norm{X}_{F}^{2} \nonumber \\ 
&\le 1.1(n+1)n\uu\norm{X}_{F}^{2}, \nonumber
\end{split}
\end{equation}
with probability at least $Q(\eta, mn^{2})$. Therefore, \eqref{eq:edc} holds. Based on the results above, Lemma~\ref{lemma 4.1} holds.
\end{proof}

\begin{lemma}
\label{lemma 4.2}
For $Y^{-1}$ and $XY^{-1}$ in \eqref{eq:49}, we have
\begin{align}
\norm{Y^{-1}}_{2} &\le \frac{1}{\sqrt{(\sigma_{min}(X))^{2}+0.9s}}, \label{eq:417} \\
\norm{XY^{-1}}_{2} &\le 1.5, \label{eq:418}
\end{align}
with probability at least $Q(\eta, mn^{2})$. 
\end{lemma}
\begin{proof}
The steps to prove \eqref{eq:417} and \eqref{eq:418} follow the same approach as that in \cite{Shifted, Columns, CSparse}. Since \eqref{eq:ecc} and \eqref{eq:edc} are used in the proof, \eqref{eq:417} and \eqref{eq:418} hold with probability at least $Q(\eta, mn^{2})$. Thus, Lemma~\ref{lemma 4.2} holds.
\end{proof}

\begin{lemma}
\label{lemma 4.3}
For $E_{WY}$ in \eqref{eq:49}, we have
\begin{equation} \label{eq:ewy2}
\norm{E_{WY}}_{2} \le 1.03n\sqrt{n}\uu \cdot \norm{W}_{2}\norm{X}_{F}, 
\end{equation}
with probability at least $Q(\eta, mn^{2})$. 
\end{lemma}
\begin{proof}
With Lemma~\ref{lemma 2.6} and \eqref{eq:49}, we can have
\begin{equation} \label{eq:ewy1}
\begin{split}
\norm{E_{WY}}_{2} &\le \gamma_{n} \cdot \norm{W}_{F}\norm{Y}_{F} \\ &\le 1.02n\sqrt{n}\uu \cdot \norm{W}_{2}\norm{Y}_{F}. 
\end{split}
\end{equation}
Similar to the steps of \eqref{eq:329}, with \eqref{eq:47}, \eqref{eq:48}, \eqref{eq:44}, \eqref{eq:ecc} and \eqref{eq:edc}, we can have
\begin{equation} \label{eq:424} 
\begin{split}
\norm{Y}_{F}^{2} &\le \norm{X}_{F}^{2}+sn+n(\norm{E_{A}}_{2}+\norm{E_{B}}_{2}) \\ &\le \norm{X}_{F}^{2}+\frac{1}{100}\norm{X}_{F}^{2}+\frac{1}{1100}\norm{X}_{F}^{2} \\ &\le 1.011\norm{X}_{F}^{2},
\end{split}
\end{equation}
with probability at least $Q(\eta, mn^{2})$. Based on \eqref{eq:424}, we can have 
\begin{equation}
\norm{Y}_{F} \le 1.006\norm{X}_{F}, \label{eq:425}
\end{equation}
with probability at least $Q(\eta, mn^{2})$. Therefore, we put \eqref{eq:425} into \eqref{eq:ewy1} and we can get \eqref{eq:ewy2}. Lemma~\ref{lemma 4.3} holds.
\end{proof}

\begin{lemma}
\label{lemma 4.4}
For $W$ in \eqref{eq:49}, we have
\begin{equation} \label{eq:434b}
\norm{W}_{2} \le 1.58, 
\end{equation}
with probability at least $Q(\eta, mn^{2})$.
\end{lemma}
\begin{proof}
Based on \eqref{eq:49}, we can have
\begin{equation} \label{eq:wfw}
\begin{split}
\norm{W}_{2} &\le \norm{XY^{-1}}_{2}+\norm{E_{WY}Y^{-1}}_{2} \\ &\le \norm{XY^{-1}}_{2}+\norm{E_{WY}}_{2}\norm{Y^{-1}}_{2}.
\end{split}
\end{equation}
With \eqref{eq:32}, \eqref{eq:kxf2}, \eqref{eq:44}, \eqref{eq:417} and \eqref{eq:ewy2}, we can have
\begin{equation} \label{eq:ewyy-1}
\begin{split}
\norm{E_{WY}}_{2}\norm{Y^{-1}}_{2} &\le \frac{1.03n\sqrt{n}\uu \cdot \norm{W}_{2}\norm{X}_{F}}{\sqrt{(\sigma_{min}(X))^{2}+0.9s}} \\ &\le \frac{1.03jn\sqrt{n}\uu\norm{X}_{2}}{\sqrt{9.9j^{2}(n+1)\uu\norm{X}_{2}^{2}}} \cdot \norm{W}_{2} \\ &\le \frac{1.03}{\sqrt{9.9}} \cdot n\sqrt{\uu} \cdot \norm{W}_{2} \\ &\le 0.05\norm{W}_{2},  
\end{split}
\end{equation}
with probability at least $Q(\eta, mn^{2})$. Therefore, we put \eqref{eq:418} and \eqref{eq:ewyy-1} into \eqref{eq:wfw} and we can have
\begin{equation} \label{eq:w2v}
\norm{W}_{2} \le 1.5+0.05\norm{W}_{2}, 
\end{equation}
with probability at least $Q(\eta, mn^{2})$. With \eqref{eq:w2v}, we can have \eqref{eq:434b}. Lemma~\ref{lemma 4.4} holds.
\end{proof}

\subsection{Proof of Theorem~\ref{theorem 4.5}}
In this part, we proof Theorem~\ref{theorem 4.5} regarding $\kappa_{2}(X)$ and $\kappa_{2}(W)$.

\begin{proof}
We have already provided a bound of $\norm{W}_{2}$ in \eqref{eq:434b}. In this part, we try to give a more accurate bound of $\norm{W}_{2}$. With \eqref{eq:47}-\eqref{eq:49}, we can have
\begin{equation}
\begin{split}
W^{\top}W &= Y^{-\top}(X+E_{WY})^{\top}(X+E_{WY})Y^{-1} \nonumber \\ 
&= Y^{-\top}X^{\top}XY^{-1}+Y^{-\top}X^{\top}E_{WY}Y^{-1} \nonumber \\ 
&+ Y^{-\top}E_{WY}^{\top}XY^{-1}+Y^{-\top}E_{WY}^{\top}E_{WY}Y^{-1} \nonumber \\
&= I-Y^{-\top}(E_{A}+E_{B}+sI)Y^{-1}+(XY^{-1})^{\top}E_{WY}Y^{-1} \nonumber \\ 
&+ Y^{-\top}E_{WY}^{\top}(XY^{-1})+Y^{-\top}E_{WY}^{\top}E_{WY}Y^{-1}. \nonumber
\end{split}
\end{equation}
Therefore, we can have
\begin{equation} \label{eq:430}
\begin{split}
\norm{W^{\top}W-I}_{2} &\le \norm{Y^{-1}}_{2}^{2}(\norm{E_{A}}_{2}+\norm{E_{B}}_{2}+s)+2\norm{Y^{-1}}_{2}\norm{XY^{-1}}_{2}\norm{E_{WY}}_{2} \\ &+ \norm{Y^{-1}}_{2}^{2}\norm{E_{WY}}_{2}^{2}. 
\end{split}
\end{equation}
In the beginning, we provide an improved bound of $\norm{E_{WY}}_{2}$ in \eqref{eq:49} compared to \eqref{eq:ewy2}. For $\norm{E_{WY}}_{F}=\norm{WY-X}_{F}$, with \eqref{eq:kxf2}, \eqref{eq:ewy2} and \eqref{eq:434b}, we can have
\begin{equation} \label{eq:res}
\begin{split}
\norm{E_{WY}}_{F} &= \norm{WY-X}_{F} \\ &\le 1.03jn\sqrt{n}\uu \cdot \norm{X}_{2} \cdot 1.58 \\ &\le 1.63jn\sqrt{n}\uu\norm{X}_{2},
\end{split}
\end{equation}
with probability at least $Q(\eta, mn^{2})$. Based on \eqref{eq:44}, \eqref{eq:ecc} and \eqref{eq:edc}, we can have $\norm{E_{A}}_{2}+\norm{E_{B}}_{2} \le 0.1s$ with probability at least $Q(\eta, mn^{2})$. Therefore, with \eqref{eq:417}, we can get
\begin{equation} \label{eq:431}  
\begin{split}
\norm{Y^{-1}}_{2}^{2}(\norm{E_{A}}_{2}+\norm{E_{B}}_{2}+s) &\le \frac{1.1s}{(\sigma_{min}(X))^{2}+0.9s} \\ &\le 1.23,  
\end{split}
\end{equation}
with probability at least $Q(\eta, mn^{2})$. Based on \eqref{eq:32}, \eqref{eq:kxf2}, \eqref{eq:44}, \eqref{eq:417}, \eqref{eq:418} and \eqref{eq:res}, we can have
\begin{equation} \label{eq:432} 
\begin{split}
2\norm{Y^{-1}}_{2}\norm{XY^{-1}}_{2}\norm{E_{WY}}_{2} &\le 2 \cdot \frac{1}{\sqrt{\sigma_{min}(X)+0.9s}} \cdot 1.5 \cdot 1.63jn\sqrt{n}\uu\norm{X}_{2} \\ &\le \frac{4.89jn\sqrt{n}\uu\norm{X}_{2}}{\sqrt{9.9j^{2}(n+1)\uu\norm{X}_{2}^{2}}} \\ &\le \frac{4.89}{\sqrt{9.9}} \cdot n\sqrt{\uu} \\ &\le 0.2, 
\end{split}
\end{equation}
with probability at least $Q(\eta, mn^{2})$. With \eqref{eq:32}, \eqref{eq:kxf2}, \eqref{eq:44}, \eqref{eq:417} and \eqref{eq:res}, we can have
\begin{equation} \label{eq:433} 
\begin{split}
\norm{Y^{-1}}_{2}^{2}\norm{E_{WY}}_{2}^{2} &\le \frac{1}{\sigma_{min}(X)+0.9s} \cdot (1.63jn\sqrt{n}\uu \cdot \norm{X}_{2})^{2} \\ &\le \frac{(1.63jn\sqrt{n}\uu\norm{X}_{2})^{2}}{9.9j^{2}(n+1)\uu\norm{X}_{2}^{2}} \\ &\le \frac{1.63^{2}}{9.9} \cdot n^{2}\uu \\ &\le 0.01,
\end{split}
\end{equation}
with probability at least $Q(\eta, mn^{2})$. We put \eqref{eq:431}-\eqref{eq:433} together into \eqref{eq:430} and we can have 
\begin{equation}
\norm{W^{\top}W-I}_{2} \le 1.44, \label{eq:qqi}
\end{equation}
with probability at least $Q(\eta, mn^{2})$. With \eqref{eq:qqi}, we can easily get
\begin{equation}
\norm{W}_{2} \le 1.57, \label{eq:434}
\end{equation}
with probability at least $Q(\eta, mn^{2})$. This is an improved bound of $\norm{W}_{2}$ compared to \eqref{eq:434b}. Since we have already estimated $\norm{W}_{2}$, we still need to estimate $\sigma_{min}(W)$ in order to evaluate $\kappa_{2}(W)$. Using Lemma~\ref{lemma 2.1} and \eqref{eq:49}, we can derive 
\begin{equation} 
\sigma_{min}(W) \ge \sigma_{min}(XY^{-1})-\norm{E_{WY}Y^{-1}}_{2}. \label{eq:438}
\end{equation}
According to \eqref{eq:417} and \eqref{eq:res}, we can have
\begin{equation} \label{eq:439}
\begin{split}
\norm{E_{WY}Y^{-1}}_{2} &\le \norm{E_{WY}}_{2}\norm{Y^{-1}}_{2} \\ &\le \frac{1.63jn\sqrt{n}\uu\norm{X}_{2}}{\sqrt{(\sigma_{min}(X))^{2}+0.9s}}, 
\end{split}
\end{equation}
with probability at least $Q(\eta, mn^{2})$. Using the same method in \cite{Shifted}, we can have
\begin{equation}
\sigma_{min}(XY^{-1}) \ge \frac{\sigma_{min}(X)}{\sqrt{(\sigma_{min}(X))^{2}+s}} \cdot 0.9, \label{eq:440}
\end{equation}
with probability at least $Q(\eta, mn^{2})$. Therefore, we put \eqref{eq:439} and \eqref{eq:440} into \eqref{eq:438} and with \eqref{eq:43}, we can have
\begin{equation} \label{eq:441}
\begin{split}
\sigma_{min}(W) &\ge \frac{0.9\sigma_{min}(X)}{\sqrt{(\sigma_{min}(X))^{2}+s}}-\frac{1.63jn\sqrt{n}\uu\norm{X}_{2}}{\sqrt{(\sigma_{min}(X))^{2}+0.9s}} \\ &\ge \frac{0.9}{\sqrt{(\sigma_{min}(X))^{2}+s}} \cdot (\sigma_{min}(X)-\frac{1.63}{0.9\sqrt{0.9}} \cdot \cdot jn\sqrt{n}\uu\norm{X}_{2}) \\ &\ge \frac{\sigma_{min}(X)}{2\sqrt{(\sigma_{min}(X))^{2}+s}} \\ &= \frac{1}{2\sqrt{1+t(\kappa_{2}(X))^{2}}},
\end{split}
\end{equation}
with probability at least $Q(\eta, mn^{2})$. Here, we let $t=\frac{s}{\norm{X}_{2}^{2}}$. With \eqref{eq:434} and \eqref{eq:441}, we can have
\begin{equation}
\kappa_{2}(W) \le 2.88\sqrt{1+t(\kappa_{2}(X))^{2}}, \nonumber
\end{equation}
with probability at least $Q(\eta, mn^{2})$. \eqref{eq:436c} is proved. 

With \eqref{eq:kxf2} and \eqref{eq:44}, we can have
\begin{equation}
\begin{split}
t &=\frac{s}{\norm{X}_{2}^{2}} \\ & \ge 11j^{2} \cdot (\eta\sqrt{m}\uu+(n+1)\uu). \label{eq:ts3}
\end{split}
\end{equation}
Similar to the steps in \cite{Shifted, Columns, CSparse}, when $\kappa_{2}(X)$ is large enough, \textit{e.g.}, $\kappa_{2}(X) \ge \uu^{-\frac{1}{2}}$, with \eqref{eq:ts3}, we can have $t(\kappa_{2}(X))^{2} \ge 11j^{2} \cdot (\eta\sqrt{m}+(n+1))>>1$. Therefore, we can get 
\begin{equation}
\sqrt{1+t(\kappa_{2}(X))^{2}} \approx \sqrt{t} \cdot \kappa_{2}(X). \label{eq:approx}
\end{equation}
With \eqref{eq:436c} and \eqref{eq:approx}, it is clear to see that
\begin{equation}
\kappa_{2}(W) \le 2.88\sqrt{t} \cdot \kappa_{2}(X), \label{eq:approxk2}
\end{equation}
with probability at least $Q(\eta, mn^{2})$. Similar to  the results in \cite{Shifted, Columns, CSparse} and with \eqref{eq:deltac2}, in order to get the sufficient condition of Shifted CholeskyQR3 regarding $\kappa_{2}(X)$, we let
\begin{equation} \label{eq:442}
\begin{split}
\kappa_{2}(W) &\le 2.88\sqrt{t} \cdot \kappa_{2}(X) \\ &\le \frac{1}{8\sqrt{mn\uu+(n+1)n\uu}}.
\end{split}
\end{equation} 
When $s=11(\eta\sqrt{m}\uu+(n+1)\uu)\norm{X}_{F}^{2}$, $t=11j^{2} \cdot (\eta\sqrt{m}\uu+(n+1)\uu)$. Based on \eqref{eq:442}, we can have 
\begin{align}
\kappa_{2}(X) &\le \frac{1}{77j\sqrt{\eta\sqrt{m}\uu+(n+1)\uu} \cdot \sqrt{mn\uu+n(n+1)\uu}}, \nonumber
\end{align}
with probability at least $Q(\eta, mn^{2})$. Therefore, \eqref{eq:2s4} is proved. Theorem~\ref{theorem 4.5} is proved.
\end{proof}

\subsection{Proof of Theorem~\ref{theorem 4.6}}
Under \eqref{eq:2s4}, we proceed to prove Theorem~\ref{theorem 4.6}. 

\begin{proof}
The proof is divided into two parts, orthogonality and residual. 

\subsubsection{Orthogonality of Shifted CholeskyQR3}
In the beginning, we consider the orthogonality. For Shifted CholeskyQR3, the same as \eqref{eq:343}, for $V$ in \eqref{eq:447}, we can have
\begin{equation}
\kappa_{2}(V) \le \sqrt{\frac{69}{59}}, \label{eq:453}
\end{equation}
with probability at least $Q(\eta, mn^{2})$. Here, the same as \eqref{eq:33}, we can have \eqref{eq:443} with probability at least $Q(\eta, mn^{2})$. 

\subsubsection{Residual of Shifted CholeskyQR3}
For the residual, with \eqref{eq:49}-\eqref{eq:452}, we can have
\begin{equation}
\begin{split}
QR &= (V+E_{7})J^{-1}(JN-E_{8}) \nonumber \\ &= (V+E_{7})N-(V+E_{7})J^{-1}E_{8} \nonumber \\ &= VN+E_{7}N-QE_{8} \nonumber \\ &= 
(W+E_{3})D^{-1}(DY-E_{4})+E_{7}N-QE_{8} \nonumber \\ &= (W+E_{3})Y-(W+E_{3})D^{-1}E_{4}+E_{7}N-QE_{8} \nonumber \\ &= WY+E_{3}Y-VE_{4}+E_{7}N-QE_{8}. \nonumber
\end{split}
\end{equation}
So it is obvious that
\begin{equation} \label{eq:454}
\begin{split}
\norm{QR-X}_{F} &\le \norm{WY-X}_{F}+\norm{E_{3}}_{F}\norm{Y}_{2}+\norm{V}_{2}\norm{E_{4}}_{F} \\ &+ \norm{E_{7}}_{F}\norm{N}_{2}+\norm{Q}_{2}\norm{E_{8}}_{F}. 
\end{split}
\end{equation}

Based on the results in \cite{Shifted, Columns} and \eqref{eq:v2}, with \eqref{eq:44}, \eqref{eq:ecc} and \eqref{eq:edc}, we can bound $\norm{Y}_{2}$ as
\begin{equation}
\norm{Y}_{2} \le 1.006\norm{X}_{2}, \label{eq:455} 
\end{equation}
with probability at least $Q(\eta, mn^{2})$. For $D$ in \eqref{eq:446}, similar to \eqref{eq:330}, \eqref{eq:455} and with \eqref{eq:434}, we can estimate $\norm{D}_{2}$ and $\norm{D}_{F}$ as
\begin{equation} \label{eq:457} 
\begin{split}
\norm{D}_{2} &\le 1.005\norm{W}_{2} \\ &\le 1.58, 
\end{split}
\end{equation}
\begin{equation} \label{eq:459} 
\begin{split}
\norm{D}_{F} &\le 1.005\norm{W}_{F} \\ &\le 1.58\sqrt{n}, 
\end{split}
\end{equation}
with probability at least $Q(\eta, mn^{2})$. The same as \eqref{eq:341}, we can get
\begin{equation}
\norm{V}_{2} \le \frac{\sqrt{69}}{8}, \label{eq:458} 
\end{equation}
with probability at least $Q(\eta, mn^{2})$. For $E_{3}$ in \eqref{eq:447}, similar to the steps to get \eqref{eq:res}, with Lemma~\ref{lemma 2.6}, \eqref{eq:459} and \eqref{eq:458}, we can bound $\norm{E_{3}}_{F}$ as
\begin{equation} \label{eq:461} 
\begin{split}
\norm{E_{3}}_{F} &\le \gamma_{n} \cdot \norm{V}_{F}\norm{D}_{F} \\ &\le 1.02n\uu \cdot \frac{\sqrt{69n}}{8} \cdot 1.58\sqrt{n} \\ &\le 1.68n^{2}\uu, 
\end{split}
\end{equation}
with probability at least $Q(\eta, mn^{2})$. For $E_{4}$ in \eqref{eq:448}, using Lemma~\ref{lemma 2.6}, \eqref{eq:kxf2}, \eqref{eq:425} and \eqref{eq:459}, we can estimate $\norm{E_{4}}_{F}$ as
\begin{equation} \label{eq:462}
\begin{split}
\norm{E_{4}}_{F} &\le \gamma_{n} \cdot \norm{D}_{F}\norm{Y}_{F} \\ &\le 1.02n\uu \cdot 1.58\sqrt{n} \cdot 1.006j\norm{X}_{2} \\ &\le 1.63jn\sqrt{n}\uu\norm{X}_{2}, 
\end{split}
\end{equation}
with probability at least $Q(\eta, mn^{2})$. For $N$ in \eqref{eq:448}, based on \eqref{eq:32}, \eqref{eq:kxf2}, \eqref{eq:425}, \eqref{eq:455}-\eqref{eq:459} and \eqref{eq:462}, $\norm{N}_{2}$ and $\norm{N}_{F}$ can be bounded as
\begin{equation} \label{eq:46n2}
\begin{split}
\norm{N}_{2} &\le \norm{D}_{2}\norm{Y}_{2}+\norm{E_{4}}_{2} \\ &\le 1.58 \cdot 1.006\norm{X}_{2}+1.63jn\sqrt{n}\uu\norm{X}_{2} \\ &\le 1.62\norm{X}_{2}, 
\end{split}
\end{equation}
\begin{equation} \label{eq:464}
\begin{split}
\norm{N}_{F} &\le \norm{D}_{2}\norm{Y}_{F}+\norm{E_{4}}_{F} \\ &\le 1.58 \cdot 1.006j\norm{X}_{2}+1.63jn\sqrt{n}\uu\norm{X}_{2} \\ &\le 1.62j\norm{X}_{2},
\end{split}
\end{equation}
with probability at least $Q(\eta, mn^{2})$. For $Q$ in \eqref{eq:457}, based on \eqref{eq:k2u} and \eqref{eq:443}, we can have
\begin{equation}
\norm{Q}_{2} \le 1.01, \label{eq:466} 
\end{equation}
with probability at least $Q(\eta, mn^{2})$. For $J$ in \eqref{eq:450}, similar to \eqref{eq:459} and with \eqref{eq:458}, we can bound $\norm{J}_{2}$ as
\begin{equation} \label{eq:j22}
\begin{split}
\norm{J}_{F} &\le 1.005\norm{V}_{F} \\ &\le 1.05\sqrt{n}, 
\end{split}
\end{equation}
with probability at least $Q(\eta, mn^{2})$. Therefore, for $E_{7}$ in \eqref{eq:451}, with Lemma~\ref{lemma 2.6}, \eqref{eq:466} and \eqref{eq:j22}, we can bound $\norm{E_{7}}_{F}$ as
\begin{equation} \label{eq:470} 
\begin{split}
\norm{E_{7}}_{F} &\le \gamma_{n} \cdot \norm{Q}_{F}\norm{J}_{F} \\ &\le 1.02n\uu \cdot 1.01\sqrt{n} \cdot 1.05\sqrt{n} \\ &\le 1.09n^{2}\uu,
\end{split}
\end{equation}
with probability at least $Q(\eta, mn^{2})$. For $E_{8}$ in \eqref{eq:452}, with Lemma~\ref{lemma 2.6}, \eqref{eq:464} and \eqref{eq:j22}, we can bound $\norm{E_{8}}_{F}$ as
\begin{equation} \label{eq:471}
\begin{split}
\norm{E_{8}}_{F} &\le \gamma_{n} \cdot \norm{J}_{F}\norm{N}_{F} \\ &\le 1.02n\uu \cdot 1.05\sqrt{n} \cdot 1.62j\norm{X}_{2} \\ &\le 1.74jn\sqrt{n}\uu\norm{X}_{2}, 
\end{split}
\end{equation}
with probability at least $Q(\eta, mn^{2})$. Therefore, we put \eqref{eq:res}, \eqref{eq:455}, \eqref{eq:458}-\eqref{eq:46n2}, \eqref{eq:466}, \eqref{eq:470} and \eqref{eq:471} into \eqref{eq:454} and we can have \eqref{eq:444}. Therefore, Theorem~\ref{theorem 4.6} holds.
\end{proof}
\begin{remark}
Theorem~\ref{theorem 4.5} and Theorem~\ref{theorem 4.6} are the key theoretical results of this work. It shows that under Assumption~\ref{assumption 1} together with a better way of error estimations, we can receive a probabilistic shifted item $s$ for Shifted CholeskyQR3. Such a probabilistic shifted item $s$ for Shifted CholeskyQR3 is the smallest one as far as we know in the existing works, which can greatly improve the applicability of the algorithm while maintaining numerical stability, see the comparisons in Table~\ref{tab:Comparisonr3} and Table~\ref{tab:Comparisonr4}. We can observe these advantages in Section~\ref{sec:experiments}. Moreover, according to \eqref{eq:44}, if we deal with $X$ with a relative small $m$, it is applicable for us to take $s=\min(11(\eta\sqrt{m}\uu+(n+1)\uu)\norm{X}_{F}^{2}, 11(m\uu+(n+1)\uu)\norm{X}_{F}^{2})$. This is a more thorough $s$ for the general cases. We can easily prove the effectiveness of $s=11(m\uu+(n+1)\uu)\norm{X}_{F}^{2}$.
\end{remark}

\section{Numerical experiments}
\label{sec:experiments}
In this section, we present several groups of numerical experiments. We primarily focus on the numerical experiments of Shifted CholeskyQR3 conducted with the probabilistic $s$ in this work. All the experiments are implemented using MATLAB R2022A on our laptop. The specifications of our laptop are detailed in Table~\ref{tab:C}. 

\begin{table}
\begin{center}
\caption{The specifications of our computer}
\begin{tabular}{c|c}
\hline
Item & Specification\\
\hline \hline
System & Windows 11 family(10.0, Version 22000) \\
BIOS & GBCN17WW \\
CPU & Intel(R) Core(TM) i5-10500H CPU @ 2.50GHz  -2.5 GHz \\
Number of CPUs / node & 12 \\
Memory size / node & 8 GB \\ 
Direct Version & DirectX 12 \\
\hline
\end{tabular}
\label{tab:C}
\end{center}
\end{table}

\subsection{Numerical stability of Shifted CholeskyQR3 with the probabilistic $s$}
\label{sec:ex1}
In this section, we focus on the numerical stability of Shifted CholeskyQR3 with different $s$. We take two different $s$, that is, the probabilistic $s$ used in this work and the improved $s$ in \cite{Columns}. For the input matrix $X \in \mathbb{R}^{m\times n}$, we primarily focus on the potential influence of $\kappa_{2}(X)$, $m$ and $n$. We construct $X$ using SVD, as described in \cite{Columns, Shifted}. We set
\begin{equation}
X=O \Sigma H^{T}, \nonumber
\end{equation}
where $O \in \mathbb{R}^{m\times m}$ and $H \in \mathbb{R}^{n\times n}$ are random orthogonal matrices and 
\begin{equation}
\Sigma = {\rm diag}(1, \sigma^{\frac{1}{n-1}}, \cdots, \sigma^{\frac{n-2}{n-1}}, \sigma) \in \mathbb{R}^{m\times n}. \nonumber
\end{equation}
Here, $0<\sigma<1$ is a positive constant. Therefore, we can have $\sigma_{1}(X)=\norm{X}_{2}=1$ and $\kappa_{2}(X)=\frac{1}{\sigma}$. 

To test the influence of $\kappa_{2}(X)$, we vary $\kappa_{2}(X)$ while fixing $m=1024$, $n=32$ and $\eta=8$. When varying $m$, we fix $n=128$, $\kappa_{2}(X)=10^{12}$ and $\eta=8$. When varying $n$, we fix $m=4096$, $\kappa_{2}(X)=10^{12}$ and $\eta=8$. To assess the influence of $\kappa_{2}(X)$, we compare the accuracy of Shifted CholeskyQR3 using the improved $s$ in Lemma~\ref{lemma 2.14} from \cite{Columns}. We take the probabilistic $s$ as $s=11\eta(\sqrt{m}\uu+(n+1)\uu)\norm{X}_{F}^{2}$ and the improved $s$ as $s=11\eta(mn\uu+n(n+1)\uu)\norm{X}_{g}^{2}$. The results of the numerical experiments are presented in Table~\ref{tab:sr}–Table~\ref{tab:po}.

\begin{table}
\caption{Shifted CholeskyQR3 with the probabilistic $s$}
\centering
\begin{tabular}{||c c c c c c||}
\hline
$\kappa_{2}(X)$ & $10^{8}$ & $10^{10}$ & $10^{12}$ & $10^{14}$ & $10^{15}$ \\
\hline
Orthogonality & $1.45e-15$ & $1.51e-15$ & $1.57e-15$ & $1.75e-15$ & $1.99e-15$ \\
\hline
Residual & $4.04e-16$ & $3.79e-16$ & $3.60e-16$ & $3.23e-16$ & $3.48e-16$ \\
\hline
\end{tabular}
\label{tab:sr}
\end{table}

\begin{table}
\caption{Shifted CholeskyQR3 with the improved $s$}
\centering
\begin{tabular}{||c c c c c c||}
\hline
$\kappa_{2}(X)$ & $10^{8}$ & $10^{10}$ & $10^{12}$ & $10^{14}$ & $10^{15}$ \\
\hline
Orthogonality & $1.32e-15$ & $1.41e-15$ & $1.38e-15$ & $1.51e-15$ & $-$ \\
\hline
Residual & $4.02e-16$ & $3.69e-16$ & $3.36e-16$ & $2.97e-16$ & $-$ \\
\hline
\end{tabular}
\label{tab:sc}
\end{table}

\begin{table}
\caption{Shifted CholeskyQR3 with the probabilistic $s$ under different $m$}
\centering
\begin{tabular}{||c c c c c c||}
\hline
$m$ & $256$ & $512$ & $1024$ & $2048$ & $4096$ \\
\hline
Orthogonality & $6.15e-15$ & $4.44e-15$ & $3.77e-15$ & $2.99e-15$ & $2.81e-15$ \\
\hline
Residual & $1.06e-15$ & $1.09e-15$ & $1.08e-15$ & $1.04e-15$ & $1.08e-15$ \\
\hline
\end{tabular}
\label{tab:pr}
\end{table}

\begin{table}
\caption{Shifted CholeskyQR3 with the probabilistic $s$ under different $n$}
\centering
\begin{tabular}{||c c c c c c||}
\hline
$n$ & $128$ & $256$ & $512$ & $1024$ & $2048$ \\
\hline
Orthogonality & $2.81e-15$ & $4.43e-15$ & $8.04e-15$ & $1.43e-14$ & $2.54e-14$ \\
\hline
Residual & $1.08e-15$ & $2.03e-15$ & $3.08e-15$ & $4.35e-15$ & $5.80e-15$ \\
\hline
\end{tabular}
\label{tab:po}
\end{table}

According to Table~\ref{tab:sr} and Table~\ref{tab:sc}, we find that Shifted CholeskyQR3 with the probabilistic $s$ is numerically stable in terms of both orthogonality and residual compared to the case with the improved $s$, as indicated by \eqref{eq:443}, \eqref{eq:444}, and the results in \cite{Columns}. Shifted CholeskyQR3 with the probabilistic $s$ shows better applicability for ill-conditioned matrices. This highlights the significance of the improved error analysis with the randomized model. Furthermore, Table~\ref{tab:pr} and Table~\ref{tab:po} show that both $m$ and $n$ do not influence the numerical stability of Shifted CholeskyQR3 with the probabilistic $s$.

\subsection{The improvements of the new theoretical results for CholeskyQR2}
In this part, we make a comparison of the theoretical results with different types of analysis and its real performances for CholeskyQR2. Regarding the input matrix, we focus on $X \in \mathbb{R}^{m\times n}$ based on SVD and we fix $\norm{X}_{2}=1$ and $\kappa_{2}(X)=10^{6}$. For the condition of $\kappa_{2}(X)$, we denote $\frac{1}{8\sqrt{mn\uu+n(n+1)\uu}}$ as the ‘\textit{Deterministic bound}’ and $\frac{1}{8j\sqrt{\eta\sqrt{m}\uu+(n+1)\uu}}$ as the ‘\textit{Probabilistic bound}’. We call the real performance of the algorithm ‘\textit{Real result}’. To test the influence of $m$, we fix $n=128$ and $\eta=8$. To test the influence of $n$, we fix $m=2048$ and $\eta=8$. Comparisons of the conditon of $\kappa_{2}(X)$ with different $m$ and $n$ are shown in Table~\ref{tab:bcc1} and Table~\ref{tab:bcc2}. Similarly, regarding the residual, we denote $5n^{2}\sqrt{n}\uu\norm{X}_{2}$ as the ‘\textit{Deterministic bound}’ and $(2.30j+1.21\sqrt{n}) \cdot n\sqrt{n}\uu\norm{X}_{2}$ as the ‘\textit{Probabilistic bound}’. The real residual of CholeskyQR2 is named as ‘\textit{Real error}’. We vary $m$ and $n$ and comparisons of the residual are shown in Table~\ref{tab:bcc3} and Table~\ref{tab:bcc4}.

\begin{table}
\caption{Comparison of the condition of $\kappa_{2}(X)$ for CholeskyQR2 when $\kappa_{2}(X)=10^{6}$, $n=128$ and $\eta=8$}
\centering
\begin{tabular}{||c c c c c c||}
\hline
$m$ & $256$ & $512$ & $1024$ & $2048$ & $4096$ \\
\hline
Real result & $\ge 10^{6}$ & $\ge 10^{6}$ & $\ge 10^{6}$ & $\ge 10^{6}$ & $\ge 10^{6}$ \\
\hline
Probabilistic bound & $4.40e+05$ & $4.00e+05$ & $3.59e+05$ & $3.18e+05$ & $2.78e+05$ \\
\hline
Deterministic bound & $5.34e+04$ & $4.14e+04$ & $3.09e+04$ & $2.25e+04$ & $1.61e+04$ \\
\hline
\end{tabular}
\label{tab:bcc1}
\end{table}

\begin{table}
\caption{Comparison of the condition of $\kappa_{2}(X)$ for CholeskyQR2 when $\kappa_{2}(X)=10^{6}$, $m=4096$ and $\eta=8$}
\centering
\begin{tabular}{||c c c c c c||}
\hline
$n$ & $128$ & $256$ & $512$ & $1024$ & $2048$ \\
\hline
Real result & $\ge 10^{6}$ & $\ge 10^{6}$ & $\ge 10^{6}$ & $\ge 10^{6}$ & $\ge 10^{6}$ \\
\hline
Probabilistic bound & $2.78e+05$ & $1.89e+05$ & $1.19e+05$ & $6.94e+04$ & $3.83e+04$ \\
\hline
Deterministic bound & $1.61e+04$ & $1.12e+04$ & $7.72e+03$ & $5.18e+03$ & $3.34e+03$ \\
\hline
\end{tabular}
\label{tab:bcc2}
\end{table}

\begin{table}
\caption{Comparison of residual for CholeskyQR2 when $\kappa_{2}(X)=10^{6}$, $n=128$ and $\eta=8$}
\centering
\begin{tabular}{||c c c c c c|}
\hline
$m$ & $256$ & $512$ & $1024$ & $2048$ & $4096$ \\
\hline
Real error & $1.41e-15$ & $1.35e-15$ & $1.35e-15$ & $1.35e-15$ & $1.34e-15$ \\
\hline
Probabilistic bound & $3.04e-12$ & $3.04e-12$ & $3.04e-12$ & $3.04e-12$ & $3.04e-12$ \\
\hline
Deterministic bound & $1.03e-10$ & $1.03e-10$ & $1.03e-10$ & $1.03e-10$ & $1.03e-10$ \\
\hline
\end{tabular}
\label{tab:bcc3}
\end{table}

\begin{table}
\caption{Comparison of residual for CholeskyQR2 when $\kappa_{2}(X)=10^{12}$, $m=4096$ and $\eta=8$}
\centering
\begin{tabular}{||c c c c c c|}
\hline
$n$ & $128$ & $256$ & $512$ & $1024$ & $2048$ \\
\hline
Real error & $1.34e-15$ & $2.60e-15$ & $3.85e-15$ & $4.99e-15$ & $6.50e-15$ \\
\hline
Probabilistic bound & $3.04e-12$ & $1.21e-11$ & $4.81e-11$ & $1.92e-10$ & $7.68e-10$ \\
\hline
Deterministic bound & $1.03e-10$ & $5.82e-10$ & $3.29e-09$ & $1.86e-08$ & $1.05e-07$ \\
\hline
\end{tabular}
\label{tab:bcc4}
\end{table}

According to Table~\ref{tab:bcc1}-Table~\ref{tab:bcc4}, we find that the theoretical results of the condition of $\kappa_{2}(X)$ and the residual have distance to the real results. However, we can find that the theoretical bounds of both the condition of $\kappa_{2}(X)$ and the residual given by the improved error analysis based on the randomized model are closer to the real results than those based on deterministic error analysis, which reflects, to some extent, the advantage of the randomized model for probabilistic error analysis. This group of numerical experiments demonstrate the effectiveness the of our improved error analysis for CholeskyQR2 with the randomized model and $\norm{\cdot}_{F}$.

\subsection{The $j$-value}
\label{sec:ex3}
In this section, we examine the $j$-value defined in \eqref{eq:kxf2}. We construct the input matrix $X$ in the same manner as described in Section~\ref{sec:ex1} based on SVD. We would like to see how large the $j$-value is compared to $\sqrt{n}$. We define $l=\frac{j}{\sqrt{n}}$, $\frac{1}{\sqrt{n}} \le l \le 1$. We test the influence of $m$ and $n$ on $l$. When varying $m$, we fix $n=128$, $\kappa_{2}(X)=10^{12}$ and $\eta=8$. For varying $n$, we fix $m=4096$, $\kappa_{2}(X)=10^{12}$ and $\eta=8$. The results of the numerical experiments are presented in Table~\ref{tab:pm}-Table~\ref{tab:pn}.

\begin{table}
\caption{$l$ with different $m$ when $n=128$, $\kappa_{2}(X)=10^{12}$ and $\eta=8$}
\centering
\begin{tabular}{||c c c c c c||}
\hline
$m$ & $256$ & $512$ & $1024$ & $2048$ & $4096$ \\
\hline
$l$ & $0.2181$ & $0.2181$ & $0.2181$ & $0.2181$ & $0.2228$ \\
\hline
\end{tabular}
\label{tab:pm}
\end{table}

\begin{table}
\caption{$l$ with different $n$ when $m=4096$, $\kappa_{2}(X)=10^{12}$ and $\eta=8$}
\centering
\begin{tabular}{||c c c c c c||}
\hline
$n$ & $128$ & $256$ & $512$ & $1024$ & $2048$ \\
\hline
$l$ & $0.1488$ & $0.1416$ & $0.1380$ & $0.1363$ & $0.1354$ \\
\hline
\end{tabular}
\label{tab:pn}
\end{table}

Table~\ref{tab:pm}–Table~\ref{tab:pn} show that the $j$-value are closely related to $n$. As $n$ increases, $l=\frac{j}{\sqrt{n}}$ decrease significantly, which aligns with the lower bound of $l$, $\frac{1}{\sqrt{n}}$. Numerical experiments demonstrate that $\norm{\cdot}_{F}$ is small enough compared to $\sqrt{n}\norm{\cdot}_{2}$ in many cases, demonstrating the advantage of taking the analysis with $\norm{\cdot}_{F}$ over the original one.

\subsection{Robustness of Shifted CholeskyQR3 with the probabilistic $s$}
\label{sec:robustness}
In this section, we show the robustness of Shifted CholeskyQR3 with the probabilistic $s$. Some of the numerical experiments, to the best of our knowledge, have not been conducted in other similar works before. We construct the input matrix $X$ in the same way as described in Section~\ref{sec:ex1} and test the potential influence of $\kappa_{2}(X)$, $m$ and $n$ which is the same as those in Section~\ref{sec:ex1}. When varying $\kappa_{2}(X)$, we fix $m=1024$, $n=32$ and $\eta=8$. When varying $m$, we fix $n=128$, $\kappa_{2}(X)=10^{12}$ and $\eta=8$. When varying $n$, we fix $m=4096$, $\kappa_{2}(X)=10^{12}$ and $\eta=8$. We use $s=11\eta \cdot (\sqrt{m}\uu+(n+1)\uu)\norm{X}_{F}^{2}$ for Shifted CholeskyQR3. We record the number of successful outcomes every $30$ trials as ‘\textit{times}’ after conducting several groups and calculate the average. Numerical results are listed in Table~\ref{tab:r1}-Table~\ref{tab:r3}. Moreover, we focus on the step of $G=X^{\top}X+E_{A}$. We denote $1.1mn\uu\norm{X}_{g}^{2}$ as the ‘\textit{Deterministic bound}’ and $1.1\eta\sqrt{m}\uu\norm{X}_{F}^{2}$ as the ‘\textit{Probabilistic bound}’. We test the times of the real error of $\norm{E_{A}}_{F}$ which exceeds the ‘\textit{Probabilistic bound}’ in every $1000$ times and record it. It is named as ‘\textit{Unsuccessful times}’. Similar numerical experiments have been done in \cite{Stochastic, New} before. We vary $\kappa_{2}(X)$, $m$ and $n$ separately and record the results in Table~\ref{tab:r4}-Table~\ref{tab:r6}.

\begin{table}
\caption{Times of success with different $\kappa_{2}(X)$ for Shifted CholeskyQR3 when $m=1024$, $n=32$ and $\eta=8$}
\centering
\begin{tabular}{||c c c c c c||}
\hline
$\kappa_{2}(X)$ & $10^{8}$ & $10^{10}$ & $10^{12}$ & $10^{14}$ & $10^{15}$ \\
\hline
Times & $30$ & $30$ & $30$ & $30$ & $30$ \\
\hline
\end{tabular}
\label{tab:r1}
\end{table}

\begin{table}
\caption{Times of success with different $m$ for Shifted CholeskyQR3}
\centering
\begin{tabular}{||c c c c c c||}
\hline
$m$ & $256$ & $512$ & $1024$ & $2048$ & $4096$ \\
\hline
Times & $30$ & $30$ & $30$ & $30$ & $30$ \\
\hline
\end{tabular}
\label{tab:r2}
\end{table}

\begin{table}
\caption{Times of success with different $n$ for Shifted CholeskyQR3}
\centering
\begin{tabular}{||c c c c c c||}
\hline
$n$ & $128$ & $256$ & $512$ & $1024$ & $2048$ \\
\hline
Times & $30$ & $30$ & $30$ & $30$ & $30$ \\
\hline
\end{tabular}
\label{tab:r3}
\end{table}

\begin{table}
\caption{Comparison of $\norm{E_{A}}_{F}$ for Shifted CholeskyQR3 when $m=1024$, $n=32$ and $\eta=8$}
\centering
\begin{tabular}{||c c c c c c|}
\hline
$\kappa_{2}(X)$ & $10^{8}$ & $10^{10}$ & $10^{12}$ & $10^{14}$ & $10^{15}$ \\
\hline
Unsuccessful times & $0$ & $0$ & $0$ & $0$ & $0$ \\
\hline
Probabilistic bound & $4.49e-14$ & $4.04e-14$ & $3.76e-14$ & $3.57e-14$ & $3.50e-14$ \\
\hline
Deterministic bound & $1.80e-13$ & $1.62e-13$ & $1.50e-13$ & $1.43e-13$ & $1.40e-13$ \\
\hline
\end{tabular}
\label{tab:r4}
\end{table}

\begin{table}
\caption{Comparison of $\norm{E_{A}}_{F}$ for Shifted CholeskyQR3 when $\kappa_{2}(X)=10^{6}$, $n=128$ and $\eta=8$}
\centering
\begin{tabular}{||c c c c c c|}
\hline
$m$ & $256$ & $512$ & $1024$ & $2048$ & $4096$ \\
\hline
Unsuccessful times & $0$ & $0$ & $0$ & $0$ & $0$ \\
\hline
Probabilistic bound & $4.43e-14$ & $6.27e-14$ & $8.86e-14$ & $1.25e-13$ & $1.77e-13$ \\
\hline
Deterministic bound & $8.86e-14$ & $1.77e-13$ & $3.54e-13$ & $7.09e-13$ & $1.42e-12$ \\
\hline
\end{tabular}
\label{tab:r5}
\end{table}

\begin{table}
\caption{Comparison of $\norm{E_{A}}_{F}$ for Shifted CholeskyQR3 when $\kappa_{2}(X)=10^{12}$, $m=4096$ and $\eta=8$}
\centering
\begin{tabular}{||c c c c c c|}
\hline
$n$ & $128$ & $256$ & $512$ & $1024$ & $2048$ \\
\hline
Unsuccessful times & $0$ & $0$ & $0$ & $0$ & $0$ \\
\hline
Probabilistic bound & $1.77e-13$ & $3.21e-13$ & $4.58e-13$ & $6.10e-13$ & $2.35e-12$ \\
\hline
Deterministic bound & $1.42e-12$ & $2.57e-12$ & $4.88e-12$ & $9.51e-12$ & $1.88e-11$ \\
\hline
\end{tabular}
\label{tab:r6}
\end{table}

Table~\ref{tab:r1}-Table~\ref{tab:r3} demonstrate that Shifted CholeskyQR3 with the probabilistic $s$ exhibits strong robustness in our numerical examples, which is crucial for the practical application of this improved algorithm. According to Table~\ref{tab:r4}-Table~\ref{tab:r6}, we find that $\norm{E_{A}}_{F}$ has a very high probability which is almost $1$ to be smaller than the probabilistic bound if we take the proper $\eta$ based on $m$ and $n$. This lays a solid foundation to the effectiveness our utilization of the randomized model.

\section{Conclusions}
\label{sec:conclusions}
In this work, we present an improved error analysis of CholeskyQR-type algorithms with the randomized model. For the input matrix $X$, we receive a better theoretical result of $\kappa_{2}(X)$ and a tighter upper bound of residual for CholeskyQR2, together with a probabilistic shifted item $s$ for Shifted CholeskyQR3. Different from other works regarding CholeskyQR, $\norm{\cdot}_{F}$ is used in the theoretical analysis, demonstrating a new method of theoretical analysis for CholeskyQR. Numerical experiments demonstrate the improvement and the effectiveness of such a new type of analysis based on the randomized model. The probabilistic $s$ based on the randomized model and $\norm{X}_{F}$ improves the applicability of Shifted CholeskyQR3 compared to the case with the existing $s$ while maintaining numerical stability. Such a probabilistic $s$ is also robust enough after numerous experiments. 

\section*{Acknowledgments}
We would like to acknowledge the support and help of Professor Zhonghua Qiao from the Hong Kong Polytechnic University. Additionally, we are grateful for the ideas shared by Dr. Qinmeng Zou from Beijing University of Posts and Telecommunications, regarding the improved analysis based on the randomized model. We also appreciate the discussions with Professor Valeria Simoncini from University of Bologna and Professor Michael Kwok-Po Ng from Hong Kong Baptist University on this topic.

\section*{Conflict of interest}
The authors declare that they have no conflict of interest.
 
\section*{Data availability}
The authors declare that all data supporting the findings of this study are available within this article.

\bibliographystyle{plain}

\bibliography{references}

\end{document}